\documentclass[hidelinks,onefignum,onetabnum]{siamart251216}

\usepackage{booktabs} 
\usepackage{siunitx} 
\usepackage{algorithm}
\usepackage{changepage}
\usepackage[caption=false]{subfig}
\usepackage{multirow}
\usepackage[percent]{overpic}
\usepackage{enumitem}

\usepackage{tabularx}

\usepackage{graphicx}
\usepackage{bm}
\usepackage{flushend}
\usepackage{xcolor} 
\definecolor{mydarkgreen}{RGB}{0,90,0} 

\definecolor{mygray}{gray}{0.8}

\usepackage{soul}
\usepackage[normalem]{ulem}
\allowdisplaybreaks
\patchcmd\newpage{\vfil}{}{}{}
\newtheorem{rem}{Remark}[section]

\usepackage{lipsum}
\usepackage{amsfonts}
\usepackage{graphicx}
\usepackage{epstopdf}
\usepackage{algorithmic}
\ifpdf
  \DeclareGraphicsExtensions{.eps,.pdf,.png,.jpg}
\else
  \DeclareGraphicsExtensions{.eps}
\fi

\newsiamremark{remark}{Remark}
\newsiamremark{hypothesis}{Hypothesis}
\crefname{hypothesis}{Hypothesis}{Hypotheses}
\newsiamthm{claim}{Claim}

\headers{TBRBF Method for Surface Advection-Diffusion Equations}{XB. Li,  L. Ling, and YZ. Sun}

\title{Trajectory-Based RBF Collocation Method via Closest-Point  Embedding for Surface Advection-Diffusion Equations\thanks{ Submitted to 
\funding{This work was supported by the General Research Fund (GRF No. 12301824, 12300922) of the Hong Kong Research Grant Council and the Guangdong and Hong Kong Universities “1+1+1” Joint Research Collaboration Scheme (project no. 2025A0505000014). }} }

\author{Xiaobin Li \thanks{Department of Mathematics, Hong Kong Baptist University, Kowloon Tong, Hong Kong ({lixiaobin@hkbu.edu.hk}).}
\and
Leevan Ling \thanks{Department of Mathematics, Hong Kong Baptist University, Kowloon Tong, Hong Kong ({lling@hkbu.edu.hk}).}
\and
Yizhong Sun \thanks{Department of Mathematics, Hong Kong Baptist University, Kowloon Tong, Hong Kong ({bill950204@126.com} $\&$ {yzsun95@hkbu.edu.hk}).}}

\usepackage{amsopn}

\begin{document}
\maketitle
\begin{abstract}
We introduce a Trajectory-Based RBF Collocation (TBRBF) method for solving surface advection-diffusion equations on smooth, compact manifolds. TBRBF decouples advection and diffusion by applying a characteristic treatment with a Kansa-type RBF collocation method for diffusion PDE, which yields an operator-split characteristic (OSC) system comprising a characteristic ODE and a diffusion PDE. We rigorously prove the equivalence between the OSC system and the original surface PDE on manifolds by embedding the latter into a narrow band domain through the closest point mapping and its constant-along-normal extension. Using an intrinsic approach, we construct a time-continuous embedded PDE with push-forward operators in each chart of the atlas and establish its equivalence with the OSC system in the narrow band. Restricting the solution back to the manifold recovers the OSC system on manifolds, ensuring that the method introduces no operator splitting error. After the surface OSC system is obtained, it admits both extrinsic and intrinsic discretizations. Extensive numerical experiments confirm the robust stability and accuracy of the proposed method.
\end{abstract}



\section{Introduction}
The transport and dispersion of scalar quantities on curved surfaces pose challenges in computational mathematics, with applications ranging across geophysical fluid dynamics, oceanography, biological systems, and materials science \cite{Bertalmio-JCP01, Charles-JCP10}. In these fields, it is crucial to extend the classical advection‐diffusion framework to accommodate the intrinsic curvature and geometry of the domain, which requires the formulation of the surface advection‐diffusion equation.

Consider a closed, compact, oriented $C^2$ manifold $\mathcal{M}$ of codimension 1 embedded in $\mathbb{R}^d$. In this work, we focus on the case $d=3$, as the case $d=2$ only requires adapting certain components of the results for $d=3$. The unit outward-pointing normal vector at any point $\mathbf{p} \in \mathcal{M}$ is denoted by $\mathbf{\hat{N}}(\mathbf{p})$, and the tangent space at that point is represented by $\mathcal{T}_{\mathcal{M}}(\mathbf{p})$.
The surface advection-diffusion equation posed on the manifold \(\mathcal{M}\) takes the form
\begin{equation}
\begin{cases}
\displaystyle
\frac{\partial}{\partial t} u + \mathbf{v}_\mathcal{M} \cdot \nabla_{\mathcal{M}} u = \varepsilon \Delta_{\mathcal{M}}  u & \mathrm{on} \quad \mathcal{M}\times[0,T], \\
\displaystyle
u(\cdot, 0) = u_{\mathcal{M},0} & \mathrm{on} \quad \mathcal{M},
\end{cases}
\label{eq:adv_diff}
\end{equation}
where \(u(\mathbf{p},t)\) is the scalar quantity, and \(\varepsilon>0\) is the diffusion coefficient. The tangential velocity field \( \mathbf{v}_\mathcal{M}: \mathcal{M}\times[0,T] \rightarrow \mathcal{T}_{\mathcal{M}} \subset \mathbb{R}^d \) defines the advection transport along the manifold.  In addition, the velocity field is assumed to have the following regularity
\begin{equation}\label{Regularity-VM}
 \mathbf{v}_{\mathcal{M}} \in C^0([0,T];C^{k+1}(\mathcal{M})).
\end{equation}
The manifold gradient operator $\nabla_{\mathcal{M}}$ and the Laplace-Beltrami operator $\Delta_{\mathcal{M}}$ in \eqref{eq:adv_diff} are defined extrinsically as
\begin{align*}
    \nabla_{\mathcal{M}} := \mathcal{P} \nabla,\quad \quad
    \Delta_{\mathcal{M}} := \nabla_{\mathcal{M}} \cdot \nabla_{\mathcal{M}} = (\mathcal{P} \nabla) \cdot (\mathcal{P} \nabla),
\end{align*}
with the orthogonal projector $\mathcal{P}=\mathcal{P}(\mathbf{p}) := \mathbf{I}-\mathbf{\hat{N}}(\mathbf{p})\mathbf{\hat{N}}(\mathbf{p})^T$, which projects vectors from $\mathbb{R}^d$ at $\mathbf{p}\in\mathcal{M}$ to the tangent space $\mathcal{T}_{\mathcal{M}}(\mathbf{p})$.

Finite element methods constitute a major class of discretizations for
PDEs on stationary and evolving surfaces. Fitted surface
finite element methods, including evolving surface finite element methods, have been developed and analyzed for a broad range of
surface PDEs
\cite{Dziuk-AN13,DziukElliott2012,Buyang-NM17}.
Eulerian unfitted formulations, including TraceFEM, use finite element
spaces defined on a background mesh and restrict their traces to the
surface
\cite{OlshanskiiReuskenXu2014,Reusken2015TraceFEM,
Xianmin-SISC17}.
CutFEM provides a closely related unfitted framework in which
stabilization controls the conditioning and approximation properties
for arbitrary cut configurations
\cite{BurmanHansboLarson2015CutFEM,
HansboLarsonZahedi2015CharacteristicCutFEM}.
For convection-dominated surface transport, SUPG/streamline-diffusion
and characteristic finite element formulations have been proposed to
improve robustness
\cite{OlshanskiiReuskenXu2014,Feng-CICP20,
HansboLarsonZahedi2015CharacteristicCutFEM,Burman-CMAME20}.
These works establish stability and error
estimates under their respective geometric and discretization
assumptions. The focus of the present work is complementary: we derive
a continuous characteristic reformulation through a closest-point embedding and push-forward operators. After restriction to
the surface, the resulting characteristic and diffusion subproblems
may in principle be discretized by fitted surface FEM, Trace/CutFEM,
or meshfree methods; here we use a Kansa-type RBF collocation method.

Closest point and related embedding methods solve surface PDEs by using points in a narrow band around the surface and by enforcing a constant-along-normal (CAN) extension \cite{Ruuth-JCP08,Marz-SINUM12,Petras-JCP16,Petras-JCP18,petras2019least}. In many time-stepping implementations, the CAN property is imposed or refreshed at discrete times through closest-point extension or interpolation. Recent work on surface advection equations \cite{Wang-JSC22,Chun-SISC23} introduced time-continuous embedding or push-forward ideas, where tangential transport is extended into the narrow band so that the CAN property is preserved continuously in time. Closest-point calculus already provides embedded representations of surface diffusion operators \cite{Marz-SINUM12}. For the time-continuous push-forward framework considered here, the diffusion part must be transformed consistently with the advection part. Accordingly, we derive not only a push-forward velocity field but also a corresponding push-forward diffusion tensor, in order to connect the embedded equation with the original surface equation and the resulting OSC formulation.

To alleviate advective time-step restrictions and improve robustness for advection-dominated problems, a widely used strategy for handling the advection term involves trajectory-based approaches, such as characteristic methods \cite{Feng-CICP20,Zhao-NHT19} or semi-Lagrangian schemes \cite{Shankar-JCP18,Ling-JSC22,Liu-CMA25}, which track the evolution of the solution along flow trajectories.
In the context of surface advection equations, Shankar and Wright proposed mesh-free semi-Lagrangian RBF methods for transport on the sphere \cite{Shankar-JCP18}. The method is formulated in Cartesian coordinates and thus does not require a closest point mapping. For more general manifolds, semi-Lagrangian and RBF-FD approaches often use projection or closest-point-type extensions \cite{Ling-JSC22, petras2019least}. They evaluate the solution at backtracked points while maintaining a CAN property. RBF-FD methods also provide flexible local approximations of surface differential operators from scattered surface nodes and normal vectors \cite{ShankarWrightKirbyFogelson2015RBFFD,Shankar-SISC20}. However, when such methods are extended from pure advection to solve advection-diffusion equations on manifolds, it becomes necessary to evaluate the Laplace–Beltrami operator  and to clarify how the diffusion part enters the trajectory-based formulation. In \cite{Liu-CMA25}, the authors employ the RBF partition of unity closest point method for advection-diffusion PDE. The closest point treatment in \cite{Liu-CMA25} is based on a fixed-time CAN extension, so interpolation or extension back to the surface is still required at each time step, as similarly discussed in \cite{Wang-JSC22}. These issues highlight the need for a time-continuous embedding approach or a push-forward diffusion tensor to ensure a consistent and robust treatment of the CAN property throughout the computational process.

\begin{figure}
    \centering
    \includegraphics[width=0.78\textwidth]{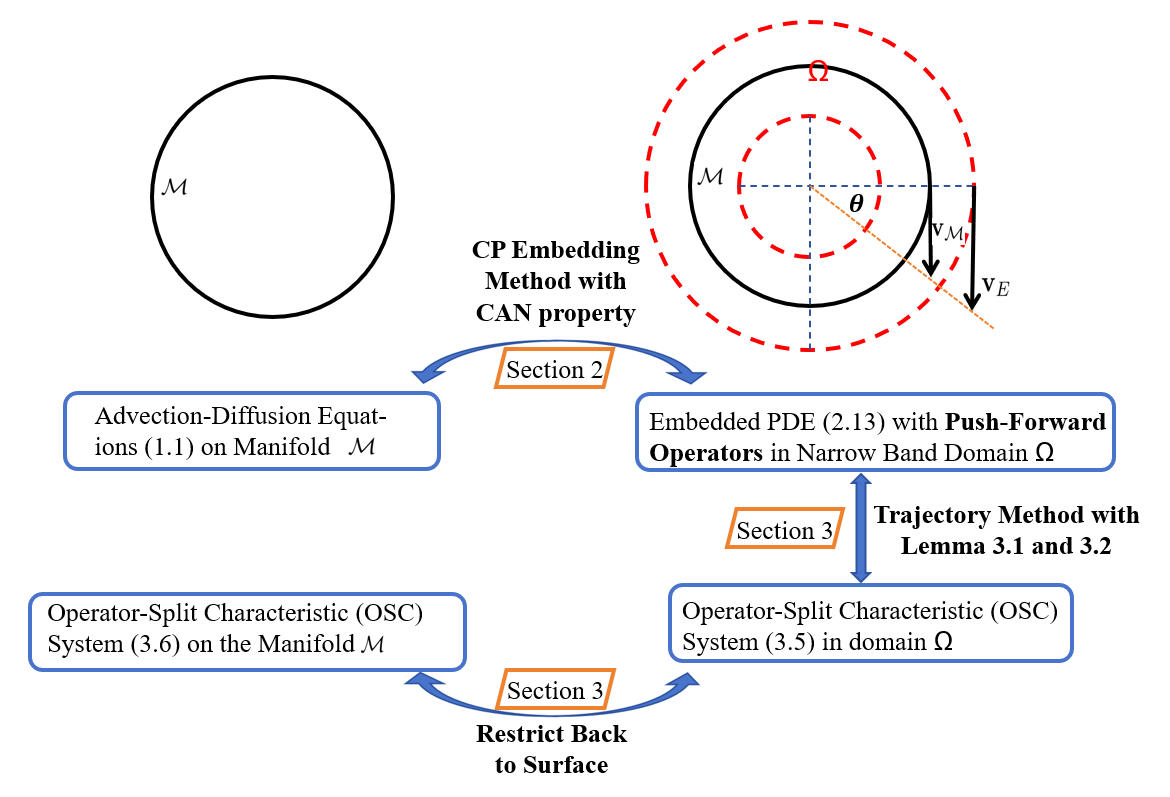}
    \caption{{Schematic diagram illustrating the equivalence between the surface PDE \eqref{eq:adv_diff} and the operator-split characteristic (OSC) system  \eqref{eq:finalTB2M}. On the left, the original surface PDE \eqref{eq:adv_diff} and the OSC system \eqref{eq:finalTB2M} are defined on the manifold $\mathcal{M}$, depicted as a circle. To establish the analytical equivalence, the manifold is embedded into the narrow band domain $\Omega$ on the right. In this embedding process, push-forward operators (e.g., velocity $
    \mathbf{v}_E$) are employed to derive the embedded PDE \eqref{eq:embedded_pde}. After proving the equivalence of the embedded PDE \eqref{eq:embedded_pde} and the OSC system \eqref{eq:final_system} in $\Omega$, restrict \eqref{eq:final_system} back to the manifold for the final result \eqref{eq:finalTB2M}.}
    }
    \label{fig:idea}
\end{figure}

As an alternative, we propose a trajectory-based RBF collocation (TBRBF) method to solve surface advection-diffusion equations. This method integrates Kansa-type RBF collocation with a trajectory-based framework, improving stability and flexibility on irregular manifolds. By directly handling the Laplace--Beltrami operator in the resulting surface OSC formulation, the TBRBF method reduces the need for repeated narrow-band interpolation associated with fixed-time CAN enforcement in classical closest point time-stepping implementations. The method avoids a separate narrow-band CAN re-extension stage; its cost is instead associated with trajectory projection, off-node RBF evaluation, and the solution of the collocation system. We refer to the coupled characteristic ODE and diffusion PDE as an operator-split characteristic (OSC) system, following the terminology used in characteristic-based discretizations. 
However, the equivalence between the OSC system on manifolds and the original surface advection-diffusion equations has not been sufficiently explored.
Previous studies \cite{Rui-NM02, ALFREDO-SINUM06} have provided analyses of the flat geometry counterpart in a bounded domain, and we aim to extend them to manifolds. 
In our work, we embed the original equation into a narrow band and rigorously establish a time-continuous embedded advection-diffusion equation with push-forward velocity and diffusion tensor in each chart of the atlas. Subsequently, we prove its equivalence with the OSC system in the narrow band domain and finally restrict back to the manifold. Note that the equivalence established above holds in the time-continuous setting, thereby eliminating any operator splitting error. Thus, a direct fixed-time closest point discretization does not by itself provide the time-continuous push-forward equivalence required in our OSC derivation, which motivates the introduction of push-forward operators. The complete schematic transformation is depicted in Figure~\ref{fig:idea}. This gives a computational framework for surface PDEs together with a continuous-level equivalence result for the trajectory form.

The rest of the paper is structured as follows.
In Section 2, we first construct an advection-aligned orthonormal frame and introduce some preliminary information about the intrinsic approach in each chart of the atlas. Then define an embedded advection-diffusion equation with a specifically designed push-forward velocity and diffusion tensor, ensuring that the embedded solution possesses the CAN property over time. Section 3 develops a trajectory-based formulation that establishes the equivalence between the embedded PDE and OSC system in the narrow band domain. We then restrict this formulation back to the surface to obtain the OSC system on the manifold. In Section 4, we detail a Kansa-type RBF collocation algorithm with various time discretization orders, and in Section 5, we present numerical experiments that demonstrate the robust stability and accuracy of the proposed TBRBF method. Finally, Section 6 concludes the paper and outlines directions for future research.

\section{The Embedded Advection-Diffusion Equations with Push-Forward Operators}
To establish the analytical equivalence between the surface PDE \eqref{eq:adv_diff} and the OSC system on the manifold $\mathcal{M}$, as illustrated in Figure \ref{fig:idea}, we begin by embedding the manifold $\mathcal{M}$ into an extended narrow band domain $\Omega$. A key requirement for this embedding is that the resulting PDE preserves the CAN property over time, following the rationale used for pure advection equations in \cite{Wang-JSC22, Chun-SISC23}. This leads us to formulate an embedded advection-diffusion equation with push-forward velocity field and diffusion tensor, ensuring that the key dynamical features of the original system are faithfully maintained within the narrow band framework.

While the push-forward velocity we derive in \eqref{pushforwardv} agrees with the result in \cite{Chun-SISC23} for pure advection, the derivation of a corresponding push-forward diffusion tensor is not addressed in that setting. Similarly, the work on curves in \cite{Wang-JSC22} provides useful motivation, but its one-dimensional geometric setting does not directly cover surface problems, especially when the advection velocity is not aligned with a principal direction. Motivated by such observations, we construct an advection-aligned orthonormal frame and adopt an intrinsic approach within each chart of the atlas that enables a systematic and continuous derivation of both the push-forward velocity and diffusion tensor.

\subsection{ Advection-aligned orthonormal frame in each chart of the atlas and intrinsic preliminaries}
We begin by constructing a well-defined orthonormal basis for the tangent space $\mathcal{T}_{\mathcal{M}}(\mathbf{p} )$ at each point $\mathbf{p} \in \mathcal{M}$, consistent with the geometric structure imposed by the surface PDE \eqref{eq:adv_diff}. 
Specifically, the tangential velocity field $\mathbf{v}_{\mathcal{M}}(\mathbf{p},t)\in \mathcal{T}_{\mathcal{M}}(\mathbf{p})$ 
can be used to define the primary tangent vector 
$\mathbf{\hat{T}}_1(\mathbf p,t)= \mathbf{v}_{\mathcal{M}}(\mathbf p,t)/\|\mathbf{v}_{\mathcal{M}}(\mathbf p,t)\|$ 
for a manifold $\mathcal{M} \subset \mathbb{R}^3$. The secondary tangent vector is then determined by 
$\mathbf{\hat{T}}_2(\mathbf p,t)= (\mathbf{\hat{N}}(\mathbf p) \times \mathbf{\hat{T}}_1(\mathbf p,t))/ \| \mathbf{\hat{N}}(\mathbf p) \times \mathbf{\hat{T}}_1(\mathbf p,t)  \|$,
with the unit normal vector $\mathbf{\hat{N}}(\mathbf p)$. 
When the time argument is clear, we suppress $(\cdot,t)$ in $\mathbf v_\mathcal M$, $\hat{\mathbf T}_1$, and $\hat{\mathbf T}_2$.

To obtain an appropriate embedded PDE, it is necessary to establish the relationship between the push-forward operators and the operators in the surface equation. 
Crucially, preserving the CAN property over time requires that the angular velocity of the embedded PDE remains the same, which in turn necessitates working within the same local parameterization.  Notably, the constancy of the angular velocity results in a discrepancy between the tangential velocity and $\mathbf{v}_{\mathcal{M}}$, while also causing variations in the diffusion speed throughout the embedding process. These considerations motivate the adoption of an intrinsic approach. 

Recall that an atlas $\mathcal{A}$ for a topological space $\mathcal{M}$ is made up of charts, each being a homeomorphism $\Phi:U\subset\mathcal{M}\to V\subset\mathbb{R}^{2}$ from an open subset $U \subset \mathcal{M}$ to an open set $V$ in the Euclidean parameter plane, typically denoted by the pair $\{ U, \Phi\}$ \cite{Lee2013ISM}. 
In our framework, we fix a point $\mathbf{p}_0\in\mathcal{M}$ and a time $t\in[0,T]$. 
We select a chart $U$ with $\mathbf{p}_{0}\in U$ and impose $\Phi(\mathbf{p}_{0})=(0,0)$.
The choice of an advection-aligned chart can depend on $t$ through $\mathbf v_\mathcal M(\cdot,t)$. 
The inverse map $\Gamma=\Phi^{-1}$ is the local surface parameterization for the chart $U$. To avoid heavy notation, we keep the symbols $\Phi$ and $\Gamma$ unchanged. All intrinsic quantities below are evaluated at the same time $t$. Within this chart we denote the coordinates by $(\theta_{1},\theta_{2})\in V$. 
We seek local coordinates $(\theta_1,\theta_2)$ whose coordinate directions are aligned with the advection-aligned orthonormal frame generated from $\mathbf v_\mathcal M(\cdot,t)$
\begin{equation}
\label{Dtheta}
\partial_{\theta_1} \Gamma(\theta_1,\theta_2) = \mathbf{\hat T}_1(\Gamma(\theta_1,\theta_2),t), \hspace{12mm} 
\partial_{\theta_2} \Gamma(\theta_1,\theta_2) = \mathbf{\hat T}_2(\Gamma(\theta_1,\theta_2),t).
\end{equation}
Condition \eqref{Dtheta} specifies the desired advection alignment, whose local existence is not automatic in general and will be formalized as part of Assumption~1 below.

\begin{figure}
    \centering
    \includegraphics[width=0.78\textwidth]{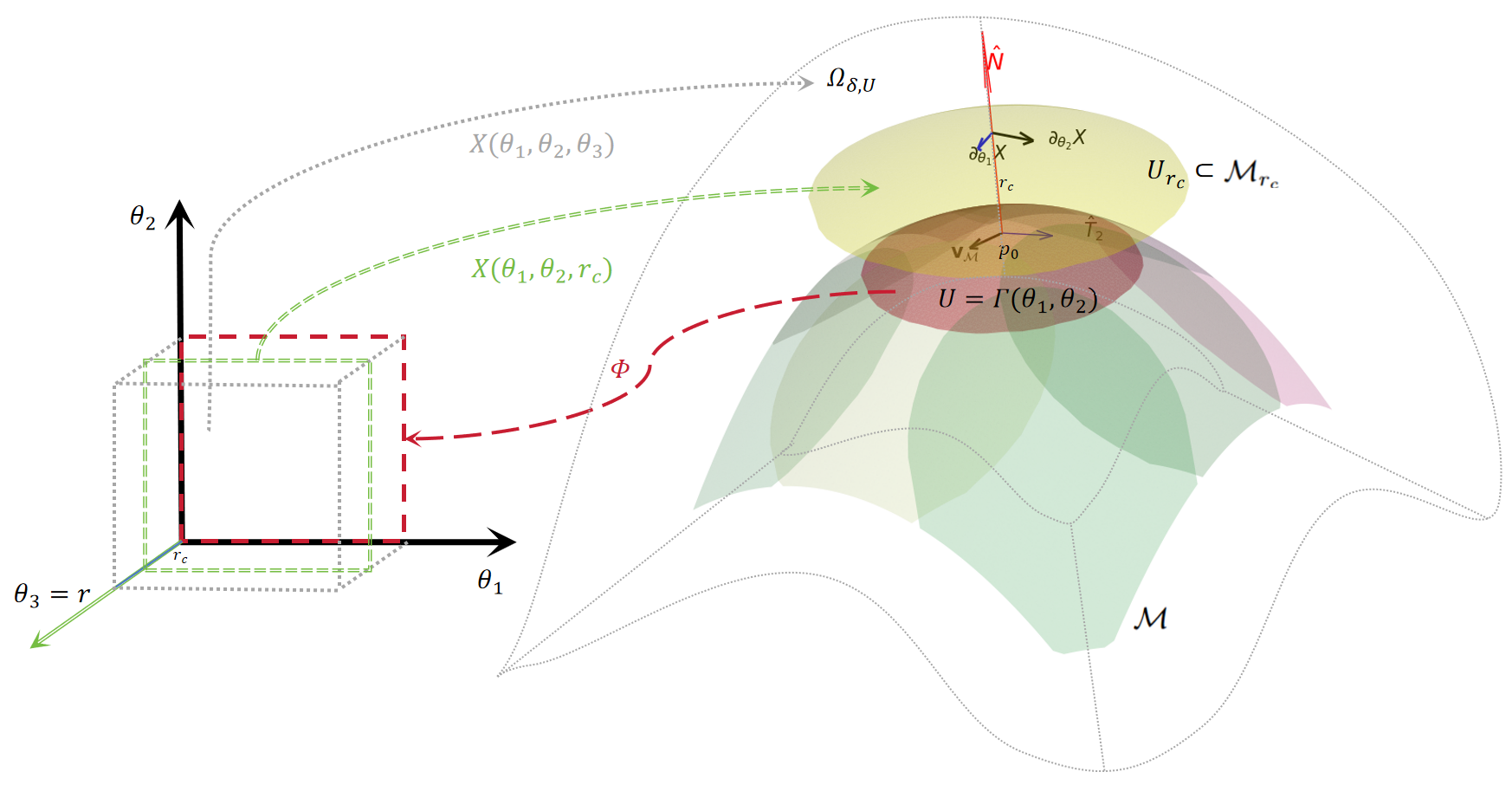}
    \caption{Construction of a local advection-aligned orthonormal frame on the surface $\mathcal{M}$ using the advection velocity as one tangent vector. The parameter domain $V$ is mapped to the surface via the atlas mapping $\Gamma=\Phi^{-1}$, where $(\theta_1,\theta_2)$ parametrize the chart $U$ on the surface $\mathcal{M}$. A third coordinate $\theta_{3}=r$ extends the chart off the surface along $\hat N$, producing the extended mapping $X$ \eqref{X-r} from $V\times[-\delta,\delta]$ (gray rectangular area) to $\Omega_{\delta,U}$ (gray dotted area), so that each level set $r=r_c:=\text{const}$ is a chart $U_{r_c}$ on the parallel surface $\mathcal M_{r_c}$.
    }
    \label{fig:pushforward}
\end{figure}

This construction incorporates the direction of the advection velocity into the coordinate system, making it especially effective for analyzing and discretizing differential operators that emerge from the surface PDE. Since this approach closely resembles a Frenet-Serret frame (TNB), which is based on curves, we refer to our local frame as the \textbf{advection-aligned orthonormal frame} $\mathbf B_2(\mathbf p,t):=\left[{\hat{\mathbf T}_{1}(\mathbf p,t),\hat{\mathbf T}_{2}(\mathbf p,t),\hat{\mathbf N}(\mathbf p)}\right] $.
When no confusion arises, we suppress the explicit dependence on $(\mathbf p,t)$. While $\mathbf B_2$ is intrinsically defined on $\mathcal M$, our construction naturally extends to the later defined narrow band domain \eqref{eq:omega_delta}. A geometric illustration is offered in \cref{fig:pushforward}.

As we will use intrinsic differential operators, we summarize a standard identity from surface differential geometry in Appendix~\ref{app:diffgeo}.
Let $L,M,N$ denote the coefficients of the second fundamental form as defined in Appendix~\ref{app:diffgeo}. 
Then to express the derivatives of the normal vector in the direction of any tangent vector along our advection-aligned coordinate directions, we use the following Weingarten equation
\begin{equation}\label{Weingarten}
 \Big[\partial_{\theta_1} \mathbf{\hat N}, \,\partial_{\theta_2} \mathbf{\hat N}\Big]  = - \Big[  \partial_{\theta_1}  \Gamma, \, \partial_{\theta_2}  \Gamma \Big] S_{\Gamma}    = - \Big[  \partial_{\theta_1}  \Gamma, \, \partial_{\theta_2}  \Gamma \Big]\begin{bmatrix} L & M \\ M & N \end{bmatrix}.
\end{equation}

Next, we embed the original surface PDE \eqref{eq:adv_diff} within a narrow band domain surrounding the manifold $\mathcal{M}$ to obtain the embedded PDE defined below \eqref{eq:embedded_pde}. The narrow band domain is defined explicitly as follows
\begin{equation}
\Omega =\Omega_{\delta} := \{ \mathbf{x} \in \mathbb{R}^d : \|\mathbf{x} - \mathbf{p}\|_{\ell_2(\mathbb{R}^d)} < \delta \text{ for some } \mathbf{p} \in \mathcal{M} \},
\label{eq:omega_delta}
\end{equation}
with respect to a positive tubular radius \( \delta \in \mathbb{R} \). Then, the closest point mapping \(\mathrm{\mathrm{cp}} : \Omega \subset \mathbb{R}^d \to \mathcal{M}\) can be well-defined in $\Omega$ (\ref{eq:omega_delta}) for a sufficiently small $\delta$ \cite{Ruuth-JCP08},  which is constructed as follows
\begin{equation}
\mathrm{\mathrm{cp}}(\mathbf{x}) = \arg \inf_{\mathbf{p} \in \mathcal{M}} \| \mathbf{p} - \mathbf{x} \|_{\ell^2(\mathbb{R}^d)}, \hspace{6mm} \mathbf{x} \in \Omega.
\end{equation}
Building on the narrow band domain $\Omega$ (\ref{eq:omega_delta}), define the signed distance function $r: \Omega\times\mathcal{M} \to \mathbb{R}$ by
\begin{equation}\label{rdist}
r = \operatorname{dist}(\mathbf{x}, \operatorname{cp}(\mathbf{x})) := \operatorname{sign}\left(\mathbf{\hat N} \cdot (\mathbf{x}-\mathrm{cp}(\mathbf{x})) \right) \| \mathbf{x}-\mathrm{cp}(\mathbf{x}) \|_{\ell^2(\mathbb{R}^d)}.
\end{equation}
The parallel surfaces of $\mathcal{M}$ with fixed distance $r=r_c:=\text{const}$ 
can be directly defined
\begin{equation}
{\mathcal{M}}_{r_c} = \{ \mathbf{x} \in \Omega : \operatorname{sign}\left(\mathbf{\hat N} \cdot (\mathbf{x}-\mathbf{p}) \right) \|\mathbf{x} - \mathbf{p}\|_{\ell^2(\mathbb{R}^d)} = r_c \hspace{1mm} \text{ for } \mathbf{p} = \operatorname{cp}(\mathbf{x}) \in \mathcal{M}\}.
\label{eq:M_delta}
\end{equation}

For notational convenience, we set $\theta_3=r$ and denote $\Omega_{\delta,U} = U \times [-\delta,\delta]$. Using the local coordinate mapping $\Gamma$ (with the same time $t$ as above) together with the closest point mapping, we define the extended mapping $X: V\times [-\delta,\delta]\rightarrow \Omega_{\delta,U} \subset \Omega$ as
\begin{equation}\label{X-r}
    X(\boldsymbol{\theta}) = \Gamma(\theta_1,\theta_2) + \theta_3 \mathbf{\hat N}(\Gamma(\theta_1,\theta_2)) \hspace{6mm} \boldsymbol{\theta}=(\theta_1,\theta_2,\theta_3),
\end{equation}
which is illustrated geometrically in Figure 2. We also have the following properties
\begin{equation*}
    \operatorname{cp}(X(\boldsymbol{\theta})) = \Gamma(\theta_1,\theta_2), \hspace{6mm} \{ X(\theta_1,\theta_2,r_c): \Gamma(\theta_1,\theta_2) \in U \} \subset \mathcal{M}_{r_c}.
\end{equation*}
Note that all $\theta$-derivatives in \eqref{relationshipXR}-\eqref{metric-G} below are taken with $t$ held fixed. Any time-dependence enters only through the advection-aligned directions used to define the chart.

Using the Weingarten equation \eqref{Weingarten}, we can differentiate $X (\theta_1,\theta_2,\theta_3)$ with respect to $\theta_i ~ (i=1,2,3)$ to obtain the tangential and normal derivatives as below
\begin{align}\label{relationshipXR}
    \left[ \partial_{\theta_1} X, \partial_{\theta_2} X,  \partial_{\theta_3} X \right] \hspace{-0.5mm}=\hspace{-0.5mm}\left[ \partial_{\theta_1} \Gamma,  \partial_{\theta_2} \Gamma,  \mathbf{\hat N} \right]\begin{bmatrix}
        1\hspace{-0.5mm}-\hspace{-0.5mm}rL & \hspace{-0.5mm}-\hspace{-0.5mm}rM & 0 \\
        \hspace{-0.5mm}-\hspace{-0.5mm}rM & 1\hspace{-0.5mm}-\hspace{-0.5mm}rN & 0 \\
        0 & 0 & 1
    \end{bmatrix} \hspace{-0.5mm}:=\hspace{-0.5mm} \left[ \mathbf{\hat T}_1,  \mathbf{\hat T}_2,  \mathbf{\hat N} \right] \mathbf{R}.
\end{align}
Equation \eqref{relationshipXR} implies that the frame of $X$, denoted by $\mathbf B_3$, is neither parallel to our proposed advection-aligned orthogonal frame $\mathbf B_2$ nor necessarily orthogonal (see Figure 2), since the matrix $\mathbf{R}$ is not a diagonal matrix. This lack of parallelism shows why an additional geometric treatment is needed when extending the curve-based construction in \cite{Wang-JSC22} to surfaces, especially when the advection velocity is not aligned with a principal direction. 
According to the definition of the metric tensor, the components of metric tensor $\mathbf{ G}_X(\boldsymbol{\theta})$ associated with $X$ are given by $ g_{ij}(\boldsymbol{\theta}) = \partial_{\theta_i} X \cdot \partial_{\theta_j} X ~ (i,j =1,2,3)$.
Hence, in our constructed frame, the full metric tensor in coordinates $(\theta_1,\theta_2,\theta_3) $ can be expressed as
\begin{equation}
    \mathbf{ G}_X(\boldsymbol{\theta}) = \mathbf{R} \mathbf{R}^{T} = \begin{bmatrix}
        (1-rL)^2+(rM)^2 & -rM(2-r(L+N)) & 0 \\
        -rM(2-r(L+N)) & (1-rN)^2+(rM)^2 & 0 \\
        0 & 0 & 1
    \end{bmatrix}.\label{metric-G}
\end{equation}

In this study, we employ three mutually related coordinate frames
\begin{itemize}
    \item Ambient Cartesian frame $\mathbf B_1$ with matrix form $\left[{\mathbf e_{1},\mathbf e_{2},\mathbf e_{3}}\right]$.
    \item Advection-aligned frame $\mathbf B_2$ on a chart $U\subset\mathcal M$ with matrix form $\left[{\hat{\mathbf T}_{1},\hat{\mathbf T}_{2},\hat{\mathbf N}}\right]$.
    \item Narrow-band frame $\mathbf B_3$ in $\Omega_{\delta,U}\subset \Omega$ with matrix form
$\left[{\partial_{\theta_{1}}X,\partial_{\theta_{2}}X,\partial_{\theta_{3}}X}\right]$.
\end{itemize}
We can use above notation to rewrite the equation \eqref{relationshipXR} and \eqref{metric-G}
\begin{equation}\label{relationXRSim}
    \mathbf{B}_{3} = \mathbf{B}_{2} \mathbf{R}, \hspace{10mm} \mathbf{G}_X = \mathbf{B}_{3}^T\mathbf{B}_{3}.
\end{equation}

As mentioned above, to localize the analysis, we require an atlas $A^*$ whose coordinates are aligned with the advection field, formalized in the following assumption.

\textbf{Assumption 1.} 
Let $\mathbf{p}_{0}\in\mathcal{M}$ be any point with advection velocity $\mathbf{v}_{\mathcal{M}}(\mathbf p_{0},t)$ in \eqref{eq:adv_diff}. 
Choose an atlas pair $\{U,\Phi\}$ (depend on $t$ through $\mathbf v_\mathcal M(\cdot,t)$) with $\mathbf p_0\in U$, $\Phi(\mathbf p_{0})=\mathbf{0}$, and its inverse $\Gamma=\Phi^{-1}$. 
Assume that $\Gamma$ satisfies the advection alignment \eqref{Dtheta} under the advection-aligned frame $\mathbf B_2(\cdot,t)$.  
All intrinsic quantities below are evaluated at the same time $t$. The chart $U$ extends to a narrow band $\Omega_{\delta,U}\supset U$ by $X$ in \eqref{X-r}. For each  $r_{c}\in(-\delta,\delta)$,  set $\mathbf x_{0}=X(0,0,r_{c})$  so that
$\operatorname{cp}(\mathbf x_{0})=\mathbf p_{0}$. Moreover, all operators admit equivalent representations in any frames via  \eqref{relationXRSim}.

Building on the preceding preliminaries of the intrinsic approach, let atlas $A^*$ satisfy Assumption 1 and denote a scalar function $u_E(\mathbf{x},t)$ in $\mathbf B_1$. Then, we can transform $u_E(\mathbf{x},t)$  from  $\mathbf{B}_1$ coordinate to $\mathbf B_3$ coordinate as $u_E(\mathbf{x},t)\circ X := \tilde u_E(\boldsymbol{\theta},t)$. Hence, for the local metric tensor $\mathbf{ G}_X$ in \eqref{metric-G}, the gradient operator can be transformed from $\mathbf{B}_1$ to $\mathbf B_3$ as \cite[Eq. (2.2)]{Dziuk-AN13}
    \begin{align}
     \nabla u_E(\mathbf{x},t) =   \sum_{i,j=1}^3  g^{ij}  \partial_{\theta_j} \tilde u_E(\boldsymbol{\theta},t) \partial_{\theta_i}X = \mathbf{ B}_{3} \mathbf{ G}_X^{-1}  \nabla_{\theta} \tilde u_E(\boldsymbol{\theta},t),\label{gradO}
    \end{align}
where the gradient $\nabla_\theta$ in frame $\mathbf B_3$ is defined as
    $\nabla_{\theta} := [\partial_{\theta_1}, \partial_{\theta_2}, \partial_{\theta_3}]^T$.

\subsection{Time-continuous embedded PDE with push-forward operators}
We present the following theorem to show the embedded PDE with the push-forward velocity and diffusion tensor.
\begin{theorem}
Let $\mathcal{M} \subset \mathbb{R}^3$ be a closed $C^2$, compact, oriented, codimension 1 surface,
$u: \mathcal{M}\times[0,T] \rightarrow \mathbb{R}$ be the solution of the surface advection-diffusion equation \eqref{eq:adv_diff}. At each point $\mathbf p_{0}\in\mathcal M$, choose an atlas $A^*$ with associated chart $(U,\Phi)$ satisfying Assumption 1. Then, on the tubular neighborhood $\Omega_{\delta,U} \supset U$, there exist a push-forward velocity field $\mathbf{v}_E: \Omega_{\delta,U}\times[0,T] \rightarrow \mathbb{R}^d$ and a push-forward diffusion tensor $\mathbf{A}: \Omega_{\delta,U} \times [0,T] \rightarrow \mathbb{R}^{d \times d}$ with components $a_{ij} ~ (i,j=1,2,3)$ and $a_{ij}=a_{ji} ~ (i \neq j)$, such that the constant‐along‐normal (CAN) extension 
\begin{equation} \label{eq:CAN}
    u_E(\mathbf{x}, t) = u(\mathrm{cp}(\mathbf{x}), t)\quad \forall ~ \mathbf{x} \in \Omega_{\delta,U}, ~ t\in[0,T].
\end{equation}
solves the embedded PDE
\begin{equation}
\begin{cases}
\displaystyle
\frac{\partial}{\partial t} u_{E}+ \mathbf{v}_E \cdot \nabla u_{E}= \varepsilon \nabla \cdot (\mathbf{A} \nabla u_{E}), & \mathrm{in }~ \Omega_{\delta,U}, \: t\in[0,T], \\
\displaystyle
u_E(\cdot, 0) = u(\mathrm{cp}(\cdot), 0), & \mathrm{in }~ \Omega_{\delta,U}, \\
\displaystyle
\partial_{\mathbf{\hat N}} u_E = 0 \quad \mathrm{and} \quad \partial_{\mathbf{\hat N}}^{(2)} u_E = 0, & \mathrm{on }~ \mathcal{M}, \: t\in[0,T].
\end{cases}
\label{eq:embedded_pde}
\end{equation}
The push-forward operators are given by
\begin{align}
   { \mathbf{v}_E(\mathbf{x},t)  } & = \mathbf{R} \hspace{0.5mm}\mathbf{v}_\mathcal{M}(\operatorname{cp}(\mathbf{x}),t)   \quad \forall ~ \mathbf{x} \in \Omega_{\delta,U}, ~ t\in[0,T],\label{pushforwardv} \\
     \mathbf{A}
    &= 
    \begin{bmatrix}
        (1-rL)^2+ (rM)^2 & -rM(2-r(L+N)) &  a_{31} \\
        -rM(2-r(L+N)) & (1-rN)^2+ (rM)^2 &  a_{32} \\
         -\int_0^r \left( \sum_{i=1}^2\partial_{\theta_i} a_{i1}  \right) & -\int_0^r \left( \sum_{i=1}^2 \partial_{\theta_i} a_{i2} \right) &  a_{33}
    \end{bmatrix},\label{pushforwardA}
\end{align}
where $\mathbf{R}$ is defined in \eqref{relationshipXR}, $r$ denotes the signed distance function \eqref{rdist}, $\theta_1, \theta_2$ are local parametric variables, and $L,M,N$ are related to the second fundamental form \eqref{H2stFF}. Moreover, $ a_{33}$ is an arbitrary value. All intrinsic quantities in \eqref{pushforwardA} (including $L,M,N$ and $\mathbf R$) are evaluated at the same time $t$ associated with $\mathbf p=\mathrm{cp}(\mathbf x)$.
\end{theorem}
\begin{rem}
    The push‐forward velocity given in \eqref{pushforwardv} is consistent with the result obtained in \cite{Chun-SISC23} for the pure advection equations, which provides an independent validation of our derivation. It is important to emphasize that the off-diagonal components ${a}_{31}$ and ${a}_{32}$ of the push-forward diffusion tensor in \eqref{pushforwardA} are, in general, nonzero for $r \neq 0$ and vanish only when $r = 0$. This indicates that, in the embedded PDE, material diffusion is not confined to planes parallel to the surface but also incorporates diffusion processes in the normal direction. These effects are not explicitly represented in the formulations of \cite{Chun-SISC23,Wang-JSC22}, which focus on pure advection or lower-dimensional geometric settings rather than a push-forward diffusion tensor for surface advection-diffusion equations.
    
    The enforcement of the CAN property \eqref{eq:CAN} ensures that the angular velocity of the embedded PDE is invariant over time. This fundamental requirement motivates the transformation of the governing equations \eqref{eq:adv_diff} (surface advection–diffusion) and \eqref{eq:embedded_pde} (embedded PDE) from the Cartesian coordinates into the constructed advection–aligned orthonormal frame. The core idea of our proof is that, by expressing both equations in this intrinsic frame, the coefficients can be directly matched, thereby yielding explicit expressions for the push-forward velocity and the diffusion tensor.
\end{rem}

Before proceeding to the proof of Theorem 2.1, we present the following lemma, which provides the transformation of the diffusion term $\nabla \cdot (\mathbf{A}(\mathbf{x},t) \nabla u_E(\mathbf{x},t))$ defined in the Cartesian coordinate system $\mathbf{B}_1$ into the local coordinate system $\mathbf B_3$.
\begin{lemma}
    Assume all in Theorem 2.1, let $\mathbf{ A}_X(\boldsymbol{\theta},t)$ be the diffusion tensor in the local coordinate system $\mathbf B_3$, satisfying $\mathbf{A}(\mathbf{x},t)=\mathbf{ B}_{3} \mathbf{ A}_X(\boldsymbol{\theta},t)\mathbf{ B}_{3}^{-1}$. The diffusion term $\nabla \cdot \left(\mathbf{A}(\mathbf{x},t) \nabla u_E(\mathbf{x},t)\right)$ under the Cartesian coordinates $\mathbf B_1$ can be expressed as the following intrinsic form in local system $\mathbf B_3$.
    \begin{align}\label{diffusion-In}
        \nabla \cdot \left(\mathbf{A}(\mathbf{x},t) \nabla u_E(\mathbf{x},t)\right) =&  \nabla_{\theta}^T \left( \mathbf{ B}_{3} \mathbf{ A}_X(\boldsymbol{\theta},t) \mathbf{ B}_{3}^{-1} \right) \mathbf{ B}_{3}^{-1} \mathbf{ B}_{3}^{-T}  \nabla_{\theta} \tilde u_E(\boldsymbol{\theta},t)\nonumber\\ &+ \operatorname{Tr} \left( \mathbf{ B}_{3}^{-T}  (\nabla_{\theta} \nabla_{\theta}^T) \tilde u_E(\boldsymbol{\theta},t) \mathbf{B}_{3}^{-1}  \left(\mathbf{ B}_{3} \mathbf{ A}_X(\boldsymbol{\theta},t) \mathbf{ B}_{3}^{-1}\right)\right).
    \end{align}
\end{lemma}
\begin{proof} (Lemma~2.2)
    Applying the product chain rule for matrix differentiation, the diffusion term $\nabla \cdot \left(\mathbf{A}(\mathbf{x},t) \nabla u_E(\mathbf{x},t)\right)$ can be decomposed into the contraction of the gradient of $\mathbf{A}$ with $\nabla u_E$ and the double contraction of $\mathbf{A}$ with the Hessian of $u_E$, i.e.
       $\nabla \cdot \left( \mathbf{A}(\mathbf{x},t) \nabla u_E(\mathbf{x},t) \right) = \nabla \left(\mathbf{A}(\mathbf{x},t) \right)  \cdot \nabla u_E(\mathbf{x},t)  + \mathbf{A}(\mathbf{x},t) :  (\nabla \nabla^T) u_E(\mathbf{x},t)$.
    Then, using the gradient operator \eqref{gradO} and relation $\mathbf{G}_X=\mathbf{B}_{3}^T \mathbf{B}_{3}$ in \eqref{relationXRSim} to transform coordinates from $\mathbf B_1$ to $\mathbf B_3$, we can easily obtain \eqref{diffusion-In}.
\end{proof}
\begin{proof}(Theorem 2.1)
    We start by recalling the CAN property \eqref{eq:CAN}: for any $ (\theta_1, \theta_2, \theta_3)$ such that $X( \theta_1, \theta_2, \theta_3):=\mathbf{x} \in \Omega_{\delta,U} \subset \Omega$ and $ (\theta_1, \theta_2)$ with $\Gamma(\theta_1, \theta_2) :=\mathbf{p} = \operatorname{cp}(\mathbf{x}) \in U \subset \mathcal{M}$, the extended solution satisfies
    $ \tilde u_E(\boldsymbol{\theta}, t) = u_E(\mathbf{x}, t) = u(\mathbf{p}, t) = \tilde u((\theta_1,\theta_2), t)$.
    This immediately implies that the extended solution $u_E$ depends only on the surface projection. In particular, the first derivatives satisfy  
    \begin{align}
        \partial_{\theta_i} \tilde u_E(\boldsymbol{\theta}, t) &= \partial_{\theta_i} \tilde u( (\theta_1,\theta_2), t), \hspace{1mm} i =1,2,\hspace{6mm}
        \partial_{\theta_3} \tilde u_E(\boldsymbol{\theta}, t) =0. \label{thetaCAN}
    \end{align}
    Furthermore, we can derive the second derivatives as
    \begin{align}
        \partial_{\theta_i\theta_j} \tilde u_E(\boldsymbol{\theta}, t) &= \partial_{\theta_i\theta_j} \tilde u( (\theta_1,\theta_2), t), \hspace{1mm} i,j =1,2, \hspace{6mm}
        \partial_{\theta_i\theta_3} \tilde u_E(\boldsymbol{\theta}, t) =0, \hspace{1mm} i=1,2,3. \label{thetaCAN2}
    \end{align}
    Subsequently, we employ the intrinsic approach to transform the embedded PDE \eqref{eq:embedded_pde} and surface advection-diffusion equation \eqref{eq:adv_diff} from the Cartesian coordinates to the advection-aligned orthonormal frame.
    Since the intrinsic reformulation below only changes spatial frames (from $\mathbf B_1$ to $\mathbf B_3$ and then to $\mathbf B_2$) at fixed physical points, the Eulerian time derivative $\partial_t u_E(\mathbf x,t)$ is unaffected by these spatial transformations. 
Using the CAN extension \eqref{eq:CAN} and noting that $\mathrm{cp}(\mathbf x)$ is time-independent (the surface $\mathcal M$ is fixed), we obtain
\begin{equation}\label{timeD}
        \frac{\partial }{\partial t} u_{E}(\mathbf{x}, t)
        =\frac{\partial}{\partial t} u(\mathrm{cp}(\mathbf{x}),t)=\frac{\partial}{\partial t} u( \mathbf{p},t).
\end{equation}
We emphasize that the intrinsic coordinates $\boldsymbol\theta$ are used only to represent spatial derivatives (with $t$ held fixed).\newline
    We then focus on transforming the advection term $\mathbf{v}_E(\mathbf{ x}, t) \cdot \nabla u_E(\mathbf{ x}, t)$ of the embedded PDE \eqref{eq:embedded_pde} into the local coordinate system. To do this, we express the velocity $\mathbf{v}_E(\mathbf{x},t)$  in the intrinsic frame as $\mathbf{v}_E  =\mathbf{ B}_{3} \left[ \mathbf{\tilde v}_E \right]_{\mathbf B_3} $, where $\left[ \cdot \right]_{\mathbf{B}_3}$ denotes the coordinate representation with respect to the basis $\mathbf{B}_3$.
    Combining with the intrinsic gradient \eqref{gradO} and equation \eqref{relationXRSim}, the original advection term transforms as
    \begin{align}\label{advXF}
         \mathbf{v}_E(\mathbf{x},t) \cdot  \nabla u_E(\mathbf{x},t)  &= \left[ \mathbf{\tilde v}_E\right]_{\mathbf B_3}  \mathbf{ B}_{3}^T \mathbf{ B}_{3} \mathbf{ G}_X^{-1}
         \nabla_\theta \tilde u_E(\boldsymbol{\theta},t)= \left[ \mathbf{\tilde v}_E^T \right]_{\mathbf B_3}  \nabla_\theta \tilde u_E(\boldsymbol{\theta},t).
    \end{align} 
    The CAN property \eqref{eq:CAN} also implies $\nabla_{\mathcal{M}} u=\nabla u$. Denote $\mathbf{v}_{\mathcal{M}} = \mathbf{B}_{2} \left[ \mathbf{\hat v}_{\mathcal{M}}\right]_{\mathbf B_2} $ with $ \left[ \mathbf{\hat v}_{\mathcal{M}}\right]_{\mathbf B_2} =\left[\|\mathbf{v}_{\mathcal{M}}\|, 0, 0\right]^T $, where the frame $\mathbf B_3$ is overlapping with $\mathbf B_2$ while $r=0$. Then, using the definition \eqref{gradO}, the advection term $\mathbf{v}_{\mathcal{M}}(\mathbf{p},t) \cdot \nabla_{\mathcal{M}} u(\mathbf{p},t) $ of the surface PDE \eqref{eq:adv_diff} can be evaluated as
    \begin{align}
        \mathbf{v}_{\mathcal{M}}(\mathbf{ p},t)  \cdot  \nabla_{\mathcal{M}} u(\mathbf{ p},t)
        &=  \left[\mathbf{\hat v}_{\mathcal{M}}^T\right]_{\mathbf B_2} \nabla_\theta  \tilde u((\theta_1,\theta_2),t).
        \label{advMF}
    \end{align} 
    As mentioned in Remark 2.1, by using the CAN property \eqref{thetaCAN} and \eqref{relationXRSim}, the advection term in the embedded PDE \eqref{eq:embedded_pde} can also be transformed into the proposed advection-aligned orthogonal frame, i.e. $\mathbf B_3:(\partial_{\theta_1}X,\partial_{\theta_2}X,\partial_{\theta_3}X) \rightarrow \mathbf B_2: (\mathbf{\hat T}_1,\mathbf{\hat T}_2,\mathbf{\hat N})$,
    \begin{align*}
        \mathbf{v}_E (\mathbf{ x}, t) \cdot \nabla u_E(\mathbf{ x}, t)
        &= \left[ \mathbf{\hat v}_E^T \right]_{\mathbf B_2} \mathbf{R}^{-1} \nabla_\theta \tilde u((\theta_1,\theta_2),t),
    \end{align*} 
    with the velocity defined in different basis $\mathbf{v}_E = \mathbf{ B}_{3} \left[\mathbf{\tilde v}_E \right]_{\mathbf B_3}  =  \mathbf{ B}_{2} \left[\mathbf{\hat v}_E\right]_{\mathbf B_2} $. Then, under the same orthogonal frame, we can match the coefficients of $\nabla_\theta u $ to get $ \left[\mathbf{\hat v}_E^T\right]_{\mathbf B_2} \mathbf{R}^{-1} =  \left[\mathbf{\hat v}_{\mathcal{M}}^T\right]_{\mathbf B_2}$.
    We can finally get $\left[\mathbf{\hat v}\right]_{\mathbf B_2}= \mathbf{R} \left[\mathbf{\hat v}_{\mathcal{M}}\right]_{\mathbf B_2}$ in our constructed frame, which is equivalent to the relation \eqref{pushforwardv} by multiplying matrix $\mathbf{B}_{2}$. \newline    
    Regarding the diffusion term, we have already derived $\nabla \cdot  (\mathbf{A}(\mathbf{ x},t) \nabla u_E(\mathbf{ x},t)) $ in the embedded domain as \eqref{diffusion-In} in Lemma 2.2.
    Similarly, with the CAN property \eqref{eq:CAN}, the diffusion term $ \Delta_{\mathcal{M}} u(\mathbf{ p},t)$ of the surface PDE \eqref{eq:adv_diff} on $\mathcal{M}$ can be derived as
    \begin{equation}
        \Delta_{\mathcal{M}} u =  \nabla \cdot (\mathbf{\tilde I}  \nabla u) = \operatorname{Tr} \left(  \mathbf{B}_{2} (\nabla_{\theta} \nabla_{\theta}^T) \tilde u  \mathbf{B}_{2}^T \mathbf{\tilde I} \right),\label{disMF}
    \end{equation}
    where the surface diffusion tensor is $\mathbf{\tilde I} = [1,0,0;~ 0,1,0;~ 0,0,\tilde a]$ with any variable $\tilde a$. For further details on the properties and derivation of $\mathbf{\tilde I}$, please refer to \cite{Chen-SINUM23}.
    Note that the diffusion term \eqref{diffusion-In} of the embedded PDE \eqref{eq:embedded_pde} remains in the frame $\mathbf B_3$. We must apply the CAN property and \eqref{thetaCAN} to transform it into the same advection-aligned frame as the diffusion term \eqref{disMF} of the surface PDE \eqref{eq:adv_diff}. Denote the diffusion tensor $\mathbf{A}(\mathbf{x},t)=\mathbf{ B}_3 \mathbf{ A}_X \mathbf{ B}_3^{-1} = \mathbf{B}_2 \mathbf{\hat A} \mathbf{ B}_2^{-1}$ and use \eqref{metric-G}, \eqref{relationXRSim} to get
    \begin{align*}\label{disEF2M}
        \nabla \cdot \left(\mathbf{A} \nabla u_E\right) =&  \nabla_{\theta}^T \left( \mathbf{B}_2\mathbf{\hat A} \mathbf{ B}_2^T  \right) \mathbf{R}^{ - 1 }\mathbf{R}^{ -   T} \nabla_{\theta} \tilde u  \\
        &+ \operatorname{Tr} \left( \mathbf{B}_2^{ - T} \mathbf{R}^{ - T} (\nabla_{\theta} \nabla_{\theta}^T) \tilde u \mathbf{R}^{ - 1}  \mathbf{B}_2^{ - 1}  (\mathbf{B}_2 \mathbf{\hat A} \mathbf{B}_2^T)  \right).
    \end{align*}
    Since all the diffusion terms are both in the advection-aligned frame, we have
    \begin{align}
        \operatorname{Tr} \left( \mathbf{B}_2^{-T} \mathbf{R}^{-T} (\nabla_{\theta} \nabla_{\theta}^T) \tilde u \mathbf{R}^{-1}  \mathbf{B}_2^{-1}  (\mathbf{B}_2 \mathbf{\hat A} \mathbf{B}_2^T)  \right) &= \operatorname{Tr} \left( \mathbf{B}_2 (\nabla_{\theta} \nabla_{\theta}^T) \tilde u \mathbf{B}_2^T \mathbf{\tilde I} \right),\label{2orderD}\\
        \nabla_{\theta}^T \left( \mathbf{B}_2\mathbf{\hat A} \mathbf{ B}_2^T  \right) \mathbf{R}^{-1 }\mathbf{R}^{-T} \nabla_{\theta} \tilde u  &= 0.\label{1orderD}
    \end{align}
    Note that \eqref{1orderD} sets the coefficient of $\nabla_\theta \tilde u$ in the intrinsic form of the embedded diffusion operator to zero. This is required to match the embedded diffusion term with the surface diffusion operator \eqref{disMF} in the same advection-aligned coordinates induced by \eqref{Dtheta} (where $G_\Gamma=\mathbf I$), since \eqref{disMF} has no tangential first-order term.
    By applying the CAN properties \eqref{thetaCAN}-\eqref{thetaCAN2}, we conclude that the last row and column of the Hessian matrix $(\nabla_{\theta} \nabla_{\theta}^T) \tilde u$ vanish. Since all the matrices in the derivation are symmetric and the last row (or equivalently, the last column) of $\mathbf{R}$ is $[0,0,1]$, straightforward matrix operations for \eqref{2orderD} yield the components $a_{ij} ~ (i,j=1,2)$ of the diffusion tensor $\mathbf{A}$ as
       $a_{ij} = (\mathbf{B}_2 \mathbf{\hat A} \mathbf{B}_2^T)_{ij} =   (  \mathbf{ B}_2 \mathbf{\tilde I} \mathbf{ B}_2^T )_{ij} = ( \mathbf{ G}_X )_{ij}$.
    Based on the equation \eqref{thetaCAN}, we can match the coefficients with
        $\nabla_{\theta}^T ( \mathbf{B}_2 \mathbf{\hat A} \mathbf{ B}_2^T  )  = [0,0,\tilde a]\mathbf{R}\mathbf{R}^T = [0,0,\tilde a]$.
    With the components $a_{ij} ~ (i,j=1,2,3)$ of $\mathbf{A} = \mathbf{B}_2 \mathbf{\hat A} \mathbf{ B}_2^T $, a simple algebraic examination shows 
    \begin{align}\label{ODEa3}
        \partial_{\theta_3} a_{3i} = - \partial_{\theta_1} a_{1i}  - \partial_{\theta_2} a_{2i}\hspace{6mm} i=1,2,
    \end{align}
    and $ a_{33}$ depends solely on $\tilde a$, which means $a_{33}$ can be selected arbitrarily. Note that the equation \eqref{ODEa3} also holds at $r=0$ with  $\mathbf{A} = \mathbf{\tilde I}$, thereby providing initial conditions
    \begin{align}\label{initialC}
        a_{3i}(\theta_3=0) =0 \hspace{6mm} i=1,2.
    \end{align}
    The equation \eqref{ODEa3} combined with initial conditions \eqref{initialC} can become an ODE system for ${a}_{31}$ and ${a}_{32}$ with respect to variable $\theta_3$, with the solution
    \begin{equation}
        a_{3i} = \int_0^r \left(- \partial_{\theta_1} a_{1i}  - \partial_{\theta_2} a_{2i}\right)   \hspace{6mm} i=1,2.
    \end{equation}
   Finally, one verifies that on the manifold $\mathcal{M}$, when $r=0$, the embedded solution recovers the surface quantities
    \begin{equation}
        \mathbf{v}_E(\mathbf{p},t)=\mathbf{v}_\mathcal{M}(\mathbf{p},t), \hspace{6mm} \mathbf{A}(\mathbf{p},t) = \mathbf{\tilde I}, \quad  \forall ~ \mathbf{p} \in \mathcal{M}, ~ t \in[0,T], \quad \mathrm{for} ~  r=0.
        \label{vandA}
    \end{equation}
    For the embedding conditions in the third line of \eqref{eq:embedded_pde}, by construction of the narrow-band map \eqref{X-r}, we have $\partial_{\theta_3}X=\hat{\mathbf N}$, and thus $\partial_{\theta_3}\tilde u_E(\boldsymbol\theta,t)=\nabla u_E(\mathbf x,t)\cdot \partial_{\theta_3}X=\partial_{\hat{\mathbf N}}u_E(\mathbf x,t)$. Then, the CAN property implies $\partial_{\theta_3}\tilde u_E= 0$, which is equivalent to $\partial_{\hat{\mathbf N}}u_E=0$ on $\mathcal M$. Similarly, differentiating once more gives $\partial_{\theta_3\theta_3}\tilde u_E=\partial_{\hat{\mathbf N}}^{(2)}u_E$, and the CAN extension yields $\partial_{\hat{\mathbf N}}^{(2)}u_E=0$ on $\mathcal M$. Similar results are shown \cite{cheung2018kernel}.
\end{proof}

\section{The Trajectory-Based Method}
In the previous section, we embedded the advection-diffusion equation \eqref{eq:adv_diff}, defined on the manifold $\mathcal{M}$, into a narrow band region $\Omega$ \eqref{eq:omega_delta} and derived an embedded PDE \eqref{eq:embedded_pde} with push-forward operators \eqref{pushforwardv}-\eqref{pushforwardA}. Building on our proposed framework for demonstrating equivalence, as illustrated in Figure \ref{fig:idea}, we now seek to further decompose the embedded PDE into an operator-split characteristic (OSC) system in $\Omega$ through a trajectory-based method. This formulation encompasses an ODE governing the characteristic curves and a diffusion PDE. In this section, we demonstrate the equivalence between the embedded PDE \eqref{eq:embedded_pde} and an OSC system to be defined in \eqref{eq:final_system} below in the narrow band domain. This OSC system \eqref{eq:final_system} can ultimately be restricted back to the manifold $\mathcal{M}$, resulting in the OSC system on $\mathcal{M}$, which shows that the original surface advection-diffusion equation \eqref{eq:adv_diff} and the OSC surface system \eqref{eq:finalTB2M} are equivalent.

It is important to note that previous works, such as \cite{Rui-NM02} and \cite{ALFREDO-SINUM06}, have established the equivalence between these characteristic methods and the original domain PDE under certain assumptions. Specifically, \cite[Prop. 1]{Rui-NM02} showed that, provided the advection velocity is sufficiently regular and vanishes on the boundary of the domain, the characteristic ODE governing the trajectories is well-defined in the discrete setting. Building on this, \cite[Prop. 3.1]{ALFREDO-SINUM06} demonstrated the differentiability of the trajectory mapping, thereby rigorously completing the equivalence argument. However, as derived in the embedded PDE \eqref{eq:embedded_pde} in Section 2, it is evident that our push-forward velocity field $\mathbf{v}_E$ \eqref{pushforwardv} is not zero at the boundary of the narrow band domain. We can only conclude that the normal component of the velocity is zero because $\mathbf{v}_E$ is tangential. This distinction necessitates a careful extension of the existing analysis to our setting (see the proof of Lemma 3.1).

We adopt the same notation for characteristic curves $ \mathcal{X}$ as in \cite{ALFREDO-SINUM06, AQbook94}. For given \( (\mathbf{x}, t) \in \Omega \times [0, T] \), the characteristic curve based on the trajectory through \( (\mathbf{x}, t) \) is the vector function \( \mathcal{X}(\cdot\,; \mathbf{x}, t) \) solving the initial value problem
\begin{equation}
\frac{\partial \mathcal{X}}{\partial \tau}(\tau; \mathbf{x}, t) = \mathbf{v}_{E}(\mathcal{X}(\tau; \mathbf{x}, t), \tau), \hspace{4mm} \mathrm{with} \, \mathrm{initial}\, \mathrm{condition}\, \mathrm{(IC):} \: \mathcal{X}(t; \mathbf{x}, t) = \mathbf{x},
\label{eq:characteristic_line}
\end{equation}
which simply specifies that for a fixed initial pair $(\mathbf{x},t)$, when the tracking time $\tau$ equals $t$, the characteristic curve starts at the point $\mathbf{x}$.
Such initial value problem \eqref{eq:characteristic_line} depicts the trajectory of a material point, initially at position \(\mathbf{x}\) and time \(t\), as it moves under the influence of the velocity field \(\mathbf{v}_{E}\) to its position $\mathcal{X}$ at time \(\tau\).
In this case, as a function on \((\tau; \mathbf{x}, t)\), it is Lipschitz continuous in \( [0, T] \times \Omega \times [0, T] \), with regularity in $\mathbf{x}$ for fixed $t, \tau$.

It is well-established that the closest point mapping $\operatorname{cp}(\mathbf{x}) $ is $C^{k}$-smooth for any point $\mathbf{x}$ within the narrow band space for any $C^{k+1}$-smooth surface \cite{Marz-SINUM12}.
The advection-aligned orthonormal frame constructed in Section 2 ensures that, for each $\mathbf{x}\in \Omega$, the corresponding frame $(\mathbf{\hat T}_1, \mathbf{\hat T}_2,\mathbf{\hat N})(\mathbf{x})$ is well-defined.
Based on the expression of push-forward velocity \eqref{pushforwardv} from Theorem 2.1 and the regularity assumption \eqref{Regularity-VM} for the tangential velocity field $\mathbf{v}_{\mathcal{M}}$,  we obtain the following regularity result for the push-forward velocity $\mathbf{v}_E$ with normal vector $\hat{\mathbf N}(\mathrm{cp}(\mathbf x))$ extended to $\Omega$ by closest-point projection
\begin{equation}
\begin{cases}
\mathbf{v}_{E} \in C^0([0,T];C^{k}(\Omega)), & \\
 \mathbf{v}_{E}(\mathbf{x},t) \cdot\mathbf{\hat N} (\operatorname{cp}(\mathbf{x}))  = 0 & \forall ~ \mathbf{x} \in \Omega, ~ t \in [0,T].
\end{cases}
\label{veproperty}
\end{equation}

Drawing inspiration from \cite[Prop 1]{Rui-NM02} and \cite[Lemma 1]{ALFREDO-SINUM06}, we present the following lemma with different velocity properties to demonstrate the effectiveness of trajectory‐based methods in $\Omega$. This indicates that solving equation (\ref{eq:characteristic_line}) is well-defined across the entire domain $\Omega$ and yields a unique solution for each initial condition $\mathbf{x}$.
\begin{lemma}
Under the property (\ref{veproperty}), for any initial point $\mathbf{x} \in \Omega$ at time $t$ and any time $\tau_1, \tau_2  \in (0,T]$, the uniqueness result of \eqref{eq:characteristic_line} gives in particular that
\begin{equation}
    \mathcal{X}(\tau_2;\mathbf{x}, t) = \mathcal{X}(\tau_2; \mathcal{X}(\tau_1; \mathbf{x}, t), \tau_1) \in \Omega.\label{one2one}
\end{equation}
Indeed, by assuming $\tau_2=t$, we deduce that the mapping $\mathcal{X}(\tau_2; \cdot, t): \Omega \rightarrow \Omega$ is one-to-one, with inverse $\mathcal{X}( t; \cdot, \tau_2)$.
\end{lemma}
\begin{figure}
    \centering
    \includegraphics[width=0.5\textwidth]{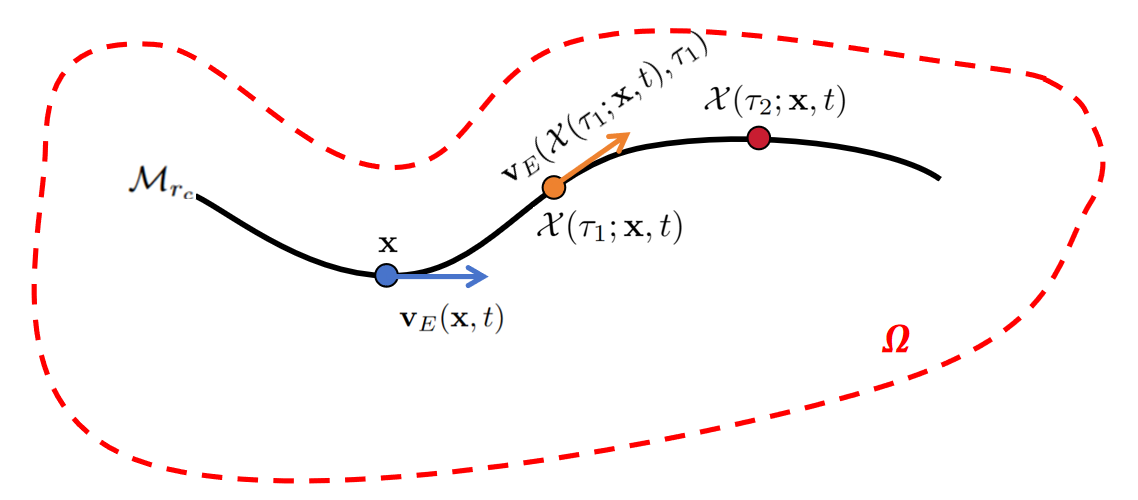}
    \caption{ Schematic illustration of the characteristic curve $\mathcal{X}(\tau;\mathbf{x}, t) $ evolving within the narrow band domain $\Omega$ (bounded by a red dashed outline). The initial point $\mathbf{x}$ (blue dot) lies on the parallel surface $\mathcal{M}_{r_c}$ (depicted as a continuous black curve). Under the tangential velocity $\mathbf{v}_E$, $\mathbf{x}$ is carried first to $\mathcal{X}(\tau_1; \mathbf{x},t)$ (orange dot) and then further to $\mathcal{X}(\tau_2; \mathcal{X}(\tau_1; \mathbf{x}, t), \tau_1) $ (red dot). This two-stage progression is equivalent to directly moving from $\mathbf{x}$ to $\mathcal{X}(\tau_2; \mathbf{x},t)$. 
    }
    \label{fig:one2one}
\end{figure}
\begin{proof}
For any initial point $\mathbf{x} \in \Omega$, including those on the boundary $\partial \Omega$, we can always consider that such a point lies within a parallel surface $\mathcal{M}_{r_c} $ \eqref{eq:M_delta}. 
There exists a smooth function $F:\Omega\rightarrow\mathbb{R}$ such that the parallel surface $\mathcal{M}_{r_c} $ is characterized by
\begin{equation*}
    \mathcal{M}_{r_c} = \{\mathbf{x}\in \Omega: ~ F(\mathbf{x}) = 0\}.
\end{equation*}
In practice, one often chooses $F$ to be the signed distance function to a reference manifold, so that $\nabla F$ is well-defined and satisfies $\|\nabla F(\mathbf{x})\|=1$ in $\Omega$.  
Consider the time derivative of \( F \) along the characteristic curve $\mathcal{X}(\tau; \mathbf{x}, t)$ defined in \eqref{eq:characteristic_line}
\begin{equation*}
    \frac{\partial}{\partial \tau} F(\mathcal{X}(\tau; \mathbf{x}, t)) =\nabla_{\mathcal{X}} F(\mathcal{X}(\tau; \mathbf{x}, t)) \cdot\frac{\partial \mathcal{X}}{\partial \tau}(\tau; \mathbf{x}, t)= \nabla_{\mathcal{X}} F(\mathcal{X}(\tau; \mathbf{x}, t)) \cdot \mathbf{v}_E(\mathcal{X}(\tau; \mathbf{x}, t), \tau),
\end{equation*}
where the operator $\nabla_{\mathcal{X}} $ denotes the gradient operator for the variable $\mathcal{X}$. Assume that the path \( \mathcal{X}(\tau; \mathbf{x}, t) \) is not exactly on the zero level set $\mathcal{M}_r$, it can lie on some other level set of $F$. At that point $\mathcal{X}$, $\nabla_{\mathcal{X}} F(\mathcal{X}(\tau; \mathbf{x}, t)) $ is the unit normal to the parallel surface (the level set of $F$) that passes through $\mathcal{X}(\tau; \mathbf{x}, t)$. 
Hence, we still have \( \mathbf{v}_E \cdot \nabla_{\mathcal{X}} F = 0 \) at $\mathcal{X}(\tau; \mathbf{x}, t)$ for all \( \tau\in [0,T] \). Since the derivative simplifies to zero, the value of $F$ is constant along the trajectory.
We now choose a specific level-set function: $F(\mathbf y)=r(\mathbf y)-r(\mathbf x)$ for $\mathbf y\in\Omega$, where $r$ is the signed distance to $\mathcal M$ in \eqref{rdist}. Then $F(\mathbf x)=0$. In $\Omega$ we have $\nabla F=\nabla r=\hat{\mathbf N}(\mathrm{cp}(\mathbf y))$. This yields $\frac{d}{d\tau}r(\mathcal X(\tau;\mathbf x,t))=\nabla r(\mathcal X)\cdot \mathbf v_E(\mathcal X,\tau)=0$. Hence $r(\mathcal X(\tau;\mathbf x,t))\equiv r(\mathbf x)$ for $\tau\in[0,T]$. Since $|r(\mathbf x)|<\delta$, we have $\mathcal X(\tau;\mathbf x,t)\in\Omega$ for all $\tau\in[0,T]$. Therefore, with the initial condition in \eqref{eq:characteristic_line}, we can get \( F(\mathcal{X}(\tau; \mathbf{x}, t)) = F(\mathcal{X}(t; \mathbf{x}, t)) = 0 \) for all \( \tau\in [0,T] \). This shows that \( \mathcal{X} \) remains on ${\mathcal{M}}_{r}\subset\Omega$ for all \( \tau\in [0,T] \) (Figure 3 provides a simple geometric explanation). 
By the uniqueness theory of ODE \cite{ODE}, the solution of \eqref{eq:characteristic_line} ensures that $\mathcal{X}(\tau_2;\mathbf{x}, t) = \mathcal{X}(\tau_2; \mathcal{X}(\tau_1; \mathbf{x}, t), \tau_1) $ for any $\tau_1,\tau_2 \in [0,T]$, which guarantees a one‐to‐one mapping.
\end{proof}

Next, based on Lemma 3.1, we refer to \cite[Prop 3.1]{ALFREDO-SINUM06} to ensure that the flow map $\mathcal{X}(\tau; \mathbf{x},t)$ is sufficiently smooth in $\mathbf{x}$, so that first and second spatial derivatives exist.
Then, using the chain rule and equation (\ref{eq:characteristic_line}), we have
\begin{align}
\frac{d }{d\tau}u_E\left(\mathcal{X}(\tau; \mathbf{x}, t), \tau\right) &= \frac{\partial u_E\left(\mathcal{X}(\tau; \mathbf{x}, t), \tau\right)}{\partial t} + \sum_{j=1}^d \frac{\partial u_E\left(\mathcal{X}(\tau; \mathbf{x}, t), \tau\right)}{\partial \mathbf{x}_j} \frac{\partial \mathcal{X}_j}{\partial \tau}\nonumber\\
&=\frac{\partial}{\partial t} u_E\left(\mathcal{X}(\tau; \mathbf{x}, t), \tau\right) +\mathbf{v}_E \cdot \nabla u_E\left(\mathcal{X}(\tau; \mathbf{x}, t), \tau\right),
\label{eq:chain_rule}
\end{align}
which is referred to as the total (material) derivative of $u_E$ with respect to $\tau\in[0,T]$ (same formulation as \cite[(12.5.4)]{AQbook94}).
By writing \eqref{eq:embedded_pde} at point \( \mathcal{X}(\tau; \mathbf{x}, t) \) and time \( \tau \), and using \eqref{eq:chain_rule}, we have the following operator-split characteristic (OSC) system
\begin{equation}
\begin{cases}
\displaystyle
\frac{\partial }{\partial \tau}\mathcal{X}(\tau; \mathbf{x}, t) = \mathbf{v}_{E}(\mathcal{X}(\tau; \mathbf{x}, t), \tau), &\mathrm{with} ~ \mathrm{IC:} \: \mathcal{X}(t; \mathbf{x}, t) = \mathbf{x},\\
\displaystyle
\frac{d u_E}{d\tau}\left(\mathcal{X}, \tau\right)  =  \varepsilon \nabla \cdot (\mathbf{A} (\mathcal{X},\tau)\nabla u_{E}(\mathcal{X},\tau)) & \text{in} \: {\Omega}, \:\tau\in[0,T], \\
\displaystyle
u_E(\cdot, 0) = u(\mathrm{cp}(\cdot), 0) & \text{in} \: \Omega,\\
\displaystyle
\partial_{\mathbf{\hat N}} u_E\left(\mathcal{X}, \tau\right) = 0 \quad \text{and} \quad \partial_{\mathbf{\hat N}}^{(2)} u_E\left(\mathcal{X}, \tau\right) = 0,
 & \text{on} \: \mathcal{M},\:\tau\in[0,T].
\end{cases}
\label{eq:final_system}
\end{equation}

Building on  \eqref{one2one} and the proof of Lemma 3.1, for any $(\mathbf x,t)\in\Omega\times[0,T]$, the trajectory $\mathcal X(\tau;\mathbf x,t)$ remains in $\Omega$ for all $\tau\in[0,T]$, hence the OSC system \eqref{eq:final_system} is well-defined in the narrow band. Moreover, the characteristic curves $\mathcal{X}$ are parallel with the manifold $\mathcal{M}$,  which implies $\operatorname{cp}(\mathcal{X}(\tau; \mathbf{x},t))= \mathcal{X}(\tau; \mathbf{p},t)$.  By restricting $\mathcal{X}(\tau; \mathbf{x},t)$ to the manifold $\mathcal{M}$ and employing the definitions of push‐forward velocity and diffusion tensor (\ref{pushforwardv})-(\ref{pushforwardA}), we can directly derive the OSC system on $\mathcal{M}$ by trajectory‐based method as
\begin{equation}
\begin{cases}\displaystyle
\frac{\partial \mathcal{X}}{\partial \tau}( \tau; \mathbf{p}, t) = \mathbf{v}_\mathcal{M}(\mathcal{X}( \tau; \mathbf{p}, t), \tau), &\mathrm{with} ~ \mathrm{IC:} \: \mathcal{X}(t;\mathbf{p}, t) = \mathbf{p}\in \mathcal{M},\\
\displaystyle
\frac{d u}{d\tau}\left(\mathcal{X}( \tau; \mathbf{p}, t), \tau\right)  =  \varepsilon \Delta_{\mathcal{M}}u(\mathcal{X}( \tau; \mathbf{p}, t),\tau)& \text{on} \: \mathcal{M}, \:\tau\in[0,T].
\end{cases}
\label{eq:finalTB2M}
\end{equation}
It is important to note that the validity of the above derivation relies on the embedded PDE with the push-forward operators preserving the CAN property of time continuity.

In summary, we have successfully established the equivalence framework illustrated in Figure \ref{fig:idea}, rigorously demonstrating the equivalence between the original surface advection-diffusion equation \eqref{eq:adv_diff} and the OSC surface system \eqref{eq:finalTB2M} on the manifold $\mathcal{M}$.
The OSC system \eqref{eq:finalTB2M} splits the surface advection-diffusion problem \eqref{eq:adv_diff} into two independent stages: tracking the characteristic curves generated by the tangential velocity field, and solving a purely diffusive equation along those trajectories.
Because the advection part is treated exactly by the flow map, the remaining diffusion step can be discretised with any convenient spatial scheme—finite differences, finite elements, or mesh-free radial basis functions. Herein, we focus on the Kansa-type RBF collocation algorithm specifically adapted to the OSC system \eqref{eq:finalTB2M}.

\section{ Kansa-Type RBF Collocation Algorithm for OSC System}
Based on the rigorous derivation presented above, we have established that the advection-diffusion equation (\ref{eq:adv_diff}) on $\mathcal{M}$ is equivalent to the trajectory-based PDEs (\ref{eq:finalTB2M}). In this section, we develop efficient numerical algorithms for such trajectory-based PDEs.

The trajectory-based method often requires determining the departure point for a given spatial location at a later time. This is accomplished by backtracing the characteristic line from any point $\mathbf{p}\in \mathcal{M}$ at time $t_{n+1}$ to its origin at time $t_n$.
Referring to \cite{Ling-JSC22}, we also employ an $s$-stage surface-restricted explicit Runge‐Kutta (RK) method to solve the first ODE in (\ref{eq:finalTB2M}), which is shown as
\begin{equation}
\begin{cases}
\mathbf{k}_i = \mathbf{v}_{\mathcal{M}} \left( \mathrm{cp} \Big( \mathbf{p} - \triangle t \sum\limits_{j=1}^{i-1} a_{ij} \mathbf{k}_j \Big), t_{n+1} - c_i \triangle t \right), & i = 1, \ldots, s, \\
\mathcal{X}( t_{n};\mathbf{p}, t_{n+1})  = \mathrm{cp} \left( \mathbf{p} - \triangle t \sum\limits_{j=1}^{s} b_j \mathbf{k}_j \right) \in \mathcal{M},
\end{cases}
\label{eq:runge_kutta}
\end{equation}
where \(\triangle t := t_{n+1} - t_{n}\), and the coefficients \( a_{ij}, b_j, \) and \( c_i \) define the RK scheme. It is important to note that we also adopt the closest point mapping in the second equation of (\ref{eq:runge_kutta}) to ensure that the new approximation point remains on the manifold $\mathcal{M}$. The convergence of the surface-restricted RK scheme is proved in \cite[Thm. 2]{Ling-JSC22}.

To compute \( u(\mathbf{p},t_{n+1}) \) from the second PDE in (\ref{eq:finalTB2M}), we have many options and introduce the time stepping scheme of the following Crank-Nicolson (CN) scheme
\begin{align*}
    &\frac{u^{n+1}(\mathcal{X}(\tau_{n+1}; \mathbf{p}, t_{n+1})) - {u}^n(\mathcal{X}(\tau_{n}; \mathbf{p}, t_{n+1}))}{\triangle \tau} \\
    &\hspace{12mm} =  \frac{\varepsilon}{2}\left( \Delta_{\mathcal{M}} u^{n+1}(\mathcal{X}(\tau_{n+1}; \mathbf{p}, t_{n+1})) +\Delta_{\mathcal{M}} {u}^n(\mathcal{X}(\tau_{n}; \mathbf{p}, t_{n+1}))\right)
\end{align*}
Note that $\mathcal{X}(\tau_{n+1}; \mathbf{p}, t_{n+1})=\mathbf{p}$ with $\tau_n=t_n$ and we only have information about fixed point $\mathbf{p}$ from the previous time step. Hence, we denote $u^{n}(\mathbf{p})$ as the solution $u(\mathbf{p},t^n)$ at the point $\mathbf{p}$ for time $t^{n}$ and $\Pi u^{n}(\mathcal{X}( t_{n}; \mathbf{p}, t_{n+1}))$ as the interpolation of the value of $u^{n}(\cdot)$. With these definitions, the CN scheme is given by
\begin{equation}
\begin{cases}
\displaystyle
\tilde{u}^{n}(\mathcal{X}(t_{n}; \mathbf{p}, t_{n+1})) := \Pi u^n(\mathcal{X}( t_{n}; \mathbf{p}, t_{n+1})), \\
\displaystyle
\frac{u^{n+1}(\mathbf{p}) \hspace{-0.5mm}-\hspace{-0.5mm} \tilde{u}^n(\mathcal{X}(t_{n}; \mathbf{p}, t_{n+1}))}{\triangle t} \hspace{-0.5mm}= \hspace{-0.5mm} \frac{\varepsilon}{2}\left( \Delta_{\mathcal{M}} u^{n+1}(\mathbf{p}) \hspace{-0.5mm}+\hspace{-0.5mm}\Delta_{\mathcal{M}} \tilde{u}^{n}(\mathcal{X}(t_{n}; \mathbf{p}, t_{n+1}))\right).
\end{cases}
\label{eq:crank_nicolson}
\end{equation}
Other time stepping methods, such as the BDF scheme, can be defined similarly. For instance, we also tested BDF methods of various time orders in our experiments, though these results are not detailed here.
Note that equation (\ref{eq:crank_nicolson}) contains a term $\Delta_{\mathcal{M}} \tilde{u}^{n}(\mathcal{X}(t_{n}; \mathbf{p}, t_{n+1}))$, which requires the utilization of information from surrounding points during the spatial discretization process. For standard grid- or mesh-based discretizations, evaluating the solution and surface differential operators at backtracked points usually requires additional interpolation, projection, or unfitted/trace finite element machinery. In RBF closest point methods, such as \cite{Liu-CMA25}, the backtracking interpolation is commonly combined with a CAN extension. This often requires extension or interpolation back to the surface at each time step. In the present algorithm, the Kansa-type RBF collocation space provides a direct interpolation mechanism at the backtracked points and is therefore convenient for the OSC formulation. 

We employ symmetric positive definite (SPD) radial basis functions (RBFs) $\Phi_m:\mathbb{R}^d \times \mathbb{R}^d \rightarrow \mathbb{R}$ that are shift-invariant and radial in nature. Following the framework established in \cite{fuselier2012scattered}, we use kernel functions of the form $\Phi_m(\|\mathbf{p}-\mathbf{z}\|)$ whose Fourier transforms exhibit algebraic decay properties. 

Let the sets of quasi-uniform trial centers be \(Z = \{\mathbf{z}_1, \ldots, \mathbf{z}_{n_Z}\} \subset \mathcal{M}\), and collocation points be \(P = \{\mathbf{p}_1, \ldots, \mathbf{p}_{n_p}\} \subset \mathcal{M}\).
The fill distance \( h \) on \(\mathcal{M}\) with respect to \(Z\) is defined as $h = h_Z:= \sup_{\xi \in \mathcal{M}} \min_{\mathbf{z}_j \in Z} \|\xi - \mathbf{z}_j\|_{\ell^2(\mathbb{R}^d)}$.
Given a set of trials $Z$ and the kernel $\Phi_m$, the finite-dimensional trial space $\mathcal{U}_Z$ is defined as $
\mathcal{U}_{Z}:=\text{span} \left\lbrace \Phi_m(\cdot, \mathbf{z}_{j}) \mid \mathbf{z}_{j} \in Z\right\rbrace$,
with the fully discrete solution
$u_{Z}^{n}(\cdot) = \sum_{j=1}^{n_Z} \lambda_{j}^{n} \Phi_m(\cdot, \mathbf{z}_{j}) = \Phi_m(\cdot, Z)\boldsymbol{\lambda}^{n}$, for $ n\geq 0$.
We define the departure point at time \(t_{n}\) traced back along the characteristic line from the collocation points $P$ at time \(t_{n+1}\) as
   $ \mathcal{X}_{n}:=\mathcal{X}(t_{n}; P, t_{n+1})$.
By substituting $\mathcal{X}_n$ and $u_{Z}^{n}(\cdot)$ into the CN scheme \eqref{eq:crank_nicolson} and the initial condition $\boldsymbol{\lambda}^0 = \Phi_m(P,Z)^{-1} u_{\mathcal{M},0}(P)$ ($\Phi_m(P, Z)$ is the collocation matrix), we obtain
\begin{align}
  \left( \Phi_m(P, Z)-
  \frac{\varepsilon}{2}\triangle t \Delta_{\mathcal{M}}\Phi_m(P, Z)\right)&\boldsymbol{\lambda}^{n+1}\nonumber\\
    &=  \left(\Phi_m(\mathcal{X}_{n}, Z) + \frac{\varepsilon}{2}\triangle t \Delta_{\mathcal{M}}\Phi_m(\mathcal{X}_{n}, Z)\right)\boldsymbol{\lambda}^{n}.
    \label{fully discrete system CN}
\end{align} 
which already incorporates the handling of the interpolation mapping $\Pi$. The extrinsic formula of the Laplace-Beltrami operator $\Delta_{\mathcal{M}}$ can be referred to \cite[Eq. (4.2)]{chen2020extrinsic}.

The trajectory-based advection-diffusion equations (\ref{eq:finalTB2M}) are fully discretized utilizing CN time discretization, in conjunction with Kansa-type RBF collocation spatial discretization.  Below, we summarize the numerical framework of the Trajectory Based RBF Collocation (TBRBF) algorithm as following Algorithm 1.

\begin{center}
\begin{minipage}{\textwidth}
\vspace{1mm}
\hrule
\vspace{1mm}
\textbf{Algorithm 1:} Trajectory Based RBF Collocation (TBRBF) Algorithm
\vspace{1mm}
\hrule
\vspace{1mm}
\begin{tabularx}{\textwidth}{@{}lX@{}}
\textbf{Step 1:} & \textbf{Identify the backtraced points:}  For each surface point $\mathbf{p}_i \in P$, solve the first ODE in \eqref{eq:finalTB2M} over the time interval $[t_{n}, t_{n+1}]$ by the surface-restricted RK \eqref{eq:runge_kutta} to identify the backtraced points $\mathcal{X}_{n}:= {\mathcal{X}}(t_{n}; P, t_{n+1})$. \\[1ex]

\textbf{Step 2:} & \textbf{Evaluate the previous solution on the backtraced points:} Use the known RBF coefficient vector $\boldsymbol\lambda^n$ and the RBF discretization $u_Z^n(\cdot)$ to evaluate the solution and the Laplace--Beltrami term at the backtraced points ${\mathcal X_n}$ and assemble the right-hand side of \eqref{fully discrete system CN}. \\[1ex]

\textbf{Step 3:} & \textbf{Solve the PDE to update the nodal solution:} Solve the fully discrete formulation \eqref{fully discrete system CN} over the time interval $[t_{n}, t_{n+1}]$ to get $\boldsymbol{\lambda}^{n+1}$, then recall RBF discretization $u_{Z}^{n}(\cdot)$ to obtain the next time level solution $u_Z^{n+1}(P)$. Other numerical discretization schemes can also be applied. \\
\end{tabularx}
\vspace{1mm}
\hrule
\end{minipage}
\end{center}

\section{Numerical Experiments}\label{sec:Numerical experiments}
In this section, we conduct numerical experiments to assess the performance of the proposed TBRBF method in solving surface advection-diffusion equations on compact, smooth manifolds.
The TBRBF method is first tested on the unit circle, showcasing optimal error convergence in time and superconvergence in space. Subsequently, the algorithm is evaluated using initial conditions with singularities on both the unit circle and the torus, demonstrating that the proposed method offers enhanced stability compared to the traditional Kansa method.
Next, we investigate the two-species Turing system with advection terms, which serves as a benchmark problem for testing numerical methods for surface PDEs. This helps to further highlight the capabilities of our TBRBF method. Finally, we extend the Turing pattern simulation to complex point clouds, providing additional evidence of the algorithm's effectiveness and wide-ranging applicability.

Conservation diagnostics provide an additional accuracy check for point-based trajectory methods.
We therefore report mass and energy errors in the convergence tests when a smooth exact solution is available.
For more complex problems, we report the discrete mass and energy at $t=0$ and $T$ to assess mass conservation and energy dissipation.
For \eqref{eq:adv_diff} (without a source term) on a closed manifold, the total mass is
$M(t)=\int_{\mathcal M} u(\mathbf p,t)\,dS$.
We also define the energy
$E(t)=\frac12\int_{\mathcal M} u(\mathbf p,t)^2\,dS$.
If the tangential velocity is solenoidal, i.e., $\nabla_{\mathcal M}\!\cdot \mathbf v_{\mathcal M}=0$, then $M(t)$ is conserved and $E(t)$ is dissipated by diffusion.
We approximate these quantities by
$M_h(t)=\sum_{i=1}^{n_P} w_i\,u_h(\mathbf p_i,t)$
and
$E_h(t)=\frac12\sum_{i=1}^{n_P} w_i\,u_h(\mathbf p_i,t)^2$,
where $w_i$ are surface quadrature weights.
When an exact solution is available, we compute the corresponding exact values $M_{\mathrm{exact}}(T)$ and $E_{\mathrm{exact}}(T)$ and report the final-time errors
$M_e=|M_h(T)-M_{\mathrm{exact}}(T)|$
and
$E_e=|E_h(T)-E_{\mathrm{exact}}(T)|$.

In all numerical tests, we use the standard Whittle-Mat\'{e}rn Sobolev kernel \cite{Matérn1986} to develop the kernel required for our RBF formulation, which is shown as
\begin{equation*}
    \Phi_m(\mathbf{p},\mathbf{z}) :=
    \frac{2^{1-m+\frac{d}{2}}}{(m-\frac{d}{2})!} \|\mathbf{p}-\mathbf{z}\|_{\ell^2}^{m-\frac{d}{2}} \mathcal{K}_{m-\frac{d}{2}}(\|\mathbf{p}-\mathbf{z}\|_{\ell^2}) \quad \quad \forall ~ \mathbf{p}, \mathbf{z} \in \mathcal{M},
\end{equation*}
where $\mathcal{K}_{\nu}$ is the modified Bessel function of the second kind of order $\nu$. Furthermore, we assume that the collocation points are always the same as the trial centers in our simulations.
For all the implicit surfaces discussed later, including the torus, sphere, Bretzel2 and CPD surfaces, we present their implicit expressions along with the corresponding tangential advection velocities in Appendix \ref{equations and velocities}. Furthermore, we have used a modified version of the \textit{distmeshsurface} algorithm \cite{distmesh} to obtain trial centers or collocation points on various surfaces.
The source code is implemented in MATLAB and can be found in the repository \cite{GitHub-xbl}.

\subsection{Example 1: Convergence tests on the circle}
In this subsection, we analyze the convergence of our numerical method for approximating the advection-diffusion equation in  \eqref{eq:adv_diff}
with the advection velocity field
$\mathbf{v}(\theta) = 
\,[-\sin(\theta),\, \cos(\theta)]^T$.
For the convergence test, we select the exact solution $u^*(\theta, t) = \exp(-\varepsilon t)\cos(\theta - t) + 1$ from which the initial conditions can be directly obtained.

We employ the surface-restricted RK method \eqref{eq:runge_kutta} with various orders $s=1,2,3$ to handle backtracking points and integrate them with the CN \eqref{eq:crank_nicolson} or BDF scheme of matching time orders for temporal discretization for the PDE of \eqref{eq:final_system}. Then, we use Mat\'{e}rn kernel with different orders $m$ in \eqref{fully discrete system CN} for spatial discretization to investigate the spatial and temporal convergence of the proposed TBRBF method.
In addition, we evaluate the impact of the diffusion coefficient $\varepsilon$ by examining some different values, a larger value (\(\varepsilon = 1\)) to a much smaller value (\(\varepsilon =  1e-5\)).

To assess the temporal convergence of the TBRBF method, we set the number of collocation points $n_P = 200$, resulting in a spatial step size (fill distance) of $h = \frac{2\pi}{n_P}$. And we fix the kernel smoothness order at $m=5$. We then select various time steps $\triangle t = 0.1, 0.05, 0.01, 0.001$ and evolve the solutions until a final time of $T=1.0$. The absolute error of the $\ell^{\infty}$ norm is used to compare the numerical solution with the exact one. As illustrated in \cref{tab:combined_convergence}, TBRBF achieves the expected convergence order under the RK2-CN scheme.
\begin{table}[!htb]
  \centering
  \caption{\textbf{(Example 1)} The temporal convergence of the RK2-CN scheme for different \(\varepsilon\) when solving \eqref{eq:adv_diff} on the unit circle ( $h = \frac{2\pi}{200}$, \(m=5\), and \(T=1.0\)).}
   \label{tab:combined_convergence}
  \begin{tabular}{ccccccc}
    \toprule
    \multirow{2}{*}{$\Delta t$}
      & \multicolumn{2}{c}{$\varepsilon=1$} & \multicolumn{2}{c}{$\varepsilon=1e-3$}
      & \multicolumn{2}{c}{$\varepsilon=1e-5$} \\
    \cmidrule(lr){2-3} \cmidrule(lr){4-5} \cmidrule(lr){6-7}
      & $\ell^\infty$ error & rate & $\ell^\infty$ error & rate
      & $\ell^\infty$ error & rate \\
    \midrule
    0.1 & 3.43e-04 & -- & 4.14e-04 & -- & 4.14e-04 & -- \\
    0.05 & 8.57e-05 & 2.00 & 1.04e-04 & 1.99 & 1.04e-04 & 1.99 \\
    0.01 & 3.43e-06 & 2.00 & 4.16e-06 & 2.00 & 4.17e-06 & 2.00 \\
    0.001 & 3.43e-08 & 2.00 & 4.16e-08 & 2.00 & 4.17e-08 & 2.00\\
    \bottomrule
  \end{tabular}
\end{table}

For the spatial convergence of TBRBF, we adopt a small time step of $\triangle t=10^{-4}$ and a total time of $T=1.0$, ensuring that the observed errors are dominated by spatial discretization.  Using the BDF3 time-stepping scheme, we investigate the effects of various fill distances $h$ and kernel smoothness orders $m$. The results, summarized in \cref{tab:combined_spatial_convergence}, demonstrate that spatial discretization achieves superconvergent behavior. In \cref{tab:combined_spatial_convergence}, $M_e$ and $E_e$ denote the final-time errors with respect to the exact mass and exact energy. 
The final-time mass and energy errors $M_e$ and $E_e$ in \cref{tab:combined_spatial_convergence} exhibit clear convergence under spatial refinement.
For more complex tests (Examples~2--3), we instead report $M_h$ and $E_h$ at $t=0$ and $T$ to assess mass conservation and energy dissipation.

\begin{table}[!htb]
\centering
\caption{\textbf{(Example 1)}  Spatial convergence of TBRBF on the unit circle at $T=1.0$ with $\Delta t=10^{-4}$ using RK3--BDF3. We report the $\ell^\infty$ error, the mass error $M_e$, and the energy error $E_e$ for different $\varepsilon$ and kernel smoothness orders $m=3,4,5$. Here $M_e$ and $E_e$ are final-time errors measured against the corresponding exact mass and energy values.} 
\label{tab:combined_spatial_convergence}
\renewcommand{\arraystretch}{1.0}
\setlength{\tabcolsep}{2.5pt}
\resizebox{\textwidth}{!}{%
\begin{tabular}{cccccccccccccccccccc}
\toprule
\multirow{2}{*}{$\varepsilon$} & \multirow{2}{*}{$h$}
& \multicolumn{6}{c}{$m=3$}
& \multicolumn{6}{c}{$m=4$}
& \multicolumn{6}{c}{$m=5$} \\
\cmidrule(lr){3-8} \cmidrule(lr){9-14} \cmidrule(lr){15-20}
&
& $\ell^\infty$  & rate & {$M_e$} & rate & {$E_e$} & rate
& $\ell^\infty$  & rate & {$M_e$} & rate & {$E_e$} & rate
& $\ell^\infty$  & rate & {$M_e$} & rate & {$E_e$} & rate \\
\midrule

\multirow{4}{*}{1}
& $\frac{2\pi}{50}$  & 3.81e-01 & --   & 1.73e+00 & --   & 1.60e+00 & --   & 7.54e-03 & --   & 3.39e-02 & --   & 3.63e-02 & --   & 1.45e-04 & --   & 6.46e-04 & --   & 6.95e-04 & --   \\
& $\frac{2\pi}{80}$  & 1.07e-01 & 2.70 & 4.86e-01 & 2.70 & 5.01e-01 & 2.47 & 7.48e-04 & 4.92 & 3.36e-03 & 4.91 & 3.61e-03 & 4.91 & 5.66e-06 & 6.90 & 2.52e-05 & 6.90 & 2.71e-05 & 6.90 \\
& $\frac{2\pi}{100}$ & 5.64e-02 & 2.89 & 2.55e-01 & 2.89 & 2.68e-01 & 2.80 & 2.47e-04 & 4.96 & 1.11e-03 & 4.96 & 1.19e-03 & 4.96 & 1.20e-06 & 6.95 & 5.35e-06 & 6.95 & 5.75e-06 & 6.95 \\
& $\frac{2\pi}{120}$ & 3.30e-02 & 2.93 & 1.50e-01 & 2.94 & 1.58e-01 & 2.89 & 9.99e-05 & 4.97 & 4.49e-04 & 4.97 & 4.82e-04 & 4.97 & 3.37e-07 & 6.97 & 1.50e-06 & 6.97 & 1.62e-06 & 6.97 \\
\midrule

\multirow{4}{*}{$10^{-5}$}
& $\frac{2\pi}{50}$  & 5.76e-04 & --   & 2.02e-05 & --   & 3.09e-05 & --   & 1.57e-05 & --   & 3.40e-07 & --   & 5.24e-07 & --   & 4.04e-07 & --   & 6.47e-09 & --   & 1.01e-08 & --   \\
& $\frac{2\pi}{80}$  & 5.76e-05 & 4.90 & 5.06e-06 & 2.95 & 7.74e-06 & 2.95 & 6.17e-07 & 6.89 & 3.36e-08 & 4.92 & 5.19e-08 & 4.92 & 6.24e-09 & 8.87 & 2.45e-10 & 6.97 & 3.85e-10 & 6.95 \\
& $\frac{2\pi}{100}$ & 1.92e-05 & 4.92 & 2.61e-06 & 2.98 & 3.98e-06 & 2.98 & 1.32e-07 & 6.93 & 1.11e-08 & 4.96 & 1.72e-08 & 4.96 & 8.55e-10 & 8.90 & 5.41e-11 & 6.76 & 8.60e-11 & 6.72 \\
& $\frac{2\pi}{120}$ & 7.83e-06 & 4.92 & 1.51e-06 & 2.98 & 2.31e-06 & 2.98 & 3.71e-08 & 6.93 & 4.49e-09 & 4.98 & 6.92e-09 & 4.98 & 1.72e-10 & 8.81 & 1.01e-11 & 9.21 & 1.34e-11 & 10.19 \\
\bottomrule
\end{tabular}%
}
\end{table}

\subsection{Example 2: Modified Shu's linear test with non-smooth initial data}
In this example, we employ a modified version of Shu's linear test \cite{jiang1996efficient}, namely a non-smooth multi-wave initial condition, to further evaluate the stability of our method in comparison to the traditional Kansa method on the unit circle.  In this modification, the original discontinuous line segments are replaced with corresponding singularity points represented by a piecewise continuous trapezoidal function. As a result, the modified multi-wave initial condition consists of four components: a smooth Gaussian function, a piecewise continuous trapezoidal function, a piecewise continuous triangular function, and a smooth elliptical function
\[
u_0(\theta) = 
\begin{cases} 
 \frac{1}{6}  G(\theta/\pi, \alpha, z + \delta) +  \frac{4}{6} G(\theta/\pi, \alpha, z) + \frac{1}{6} G(\theta/\pi, \alpha, z - \delta) , & \theta \in \textstyle \bigl[\frac{1}{5}\pi, \frac{2}{5}\pi\bigr], \\
 20\theta - 35, & \theta \in [1.75, 1.8], \\
 1, & \theta \in \textstyle \bigl[1.8, 2.5\bigr], \\
-20\theta + 51, & \theta \in [2.5,2.65],\\
 1 - \bigl|2.5\theta - 8.5\bigr|, & \theta \in [3, \textstyle 3.8], \\
 \frac{1}{6}  F(\theta/\pi, \beta, b + \delta) +  \frac{4}{6}F(\theta/\pi, \beta, b) +  \frac{1}{6}F(\theta/\pi, \beta, b - \delta) , & \theta \in \textstyle \bigl[\frac{7}{5}\pi, \frac{8}{5}\pi\bigr], \\
0, & \text{otherwise}.
\end{cases}
\]
where
\[
\begin{aligned}
G(\theta, \alpha, z) = e^{-\alpha(\theta-z)^2}, \hspace{6mm}
 F(\theta, \beta, b) = \sqrt{\max \left( 1 - \beta^2 (\theta - b)^2, 0 \right)},
\end{aligned}
\]
with \( \alpha = \log 2/(36 \delta^{2}) \), \( \beta = 10 \), \( z = 0.3 \), \( b = 1.5 \), and \( \delta = 0.005 \).
For simplicity, we use the same velocity defined in Example 1 $\mathbf{v}(\theta) = 
\,[-\sin(\theta),\, \cos(\theta)]^T$. We then apply the RK2 method in conjunction with the CN scheme for time discretization, while employing the Mat\'{e}rn kernel with a smoothness order of $m = 4$ for spatial discretization. The number of collocation points is fixed at $n_P = 500$ and the time step is set to $\triangle t= 0.01$ with a final time of $T=2\pi$.

\begin{figure}[htbp]
    \centering
    \begin{tabular}{c}
        \includegraphics[width=0.9\textwidth]{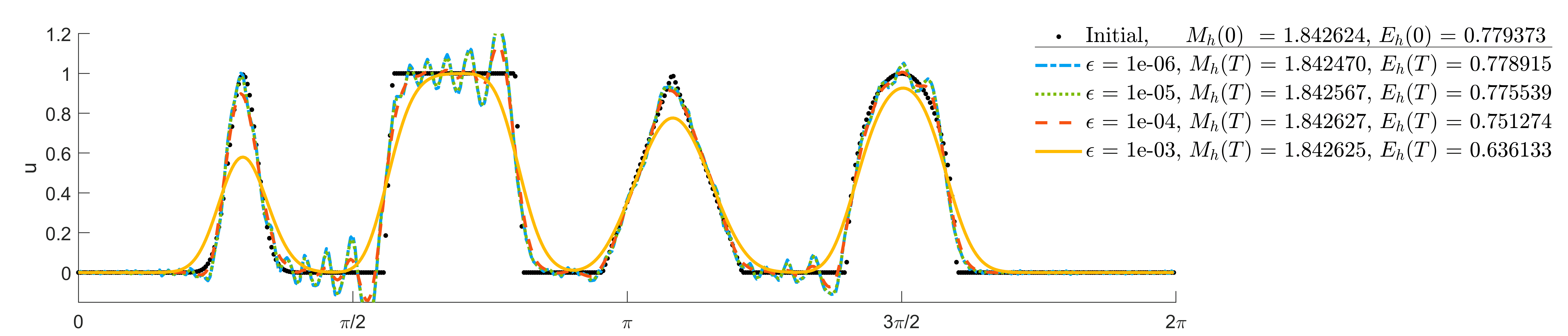}\\
        (a)  Traditional Kansa method solutions with different $\varepsilon$\\
        \vspace{0.5mm}
        \\   \includegraphics[width=0.9\textwidth]{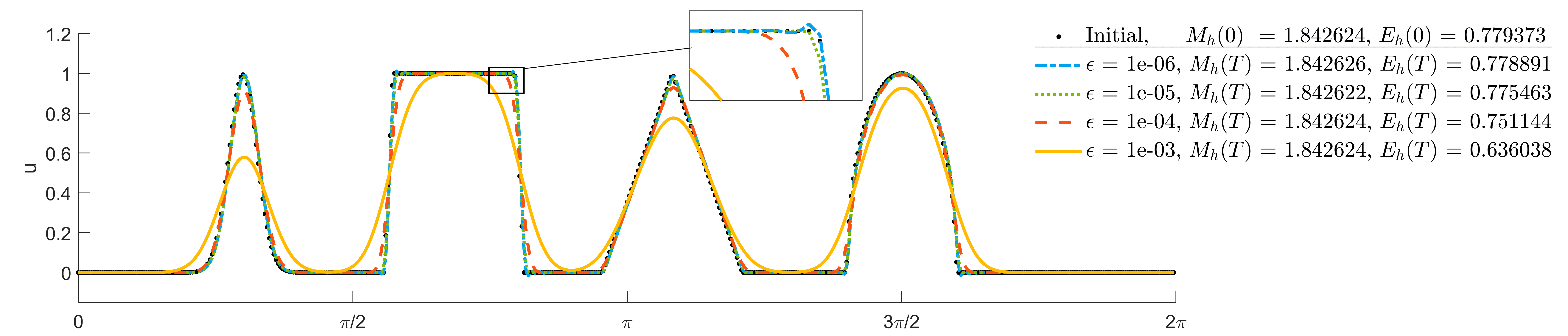}\\
          (b) TBRBF method solutions with different $\varepsilon$
    \end{tabular}
    \caption{\textbf{(Example 2)}
              Comparison between traditional Kansa method (top) and TBRBF method (bottom) to solve \eqref{eq:adv_diff} on the unit circle with non-smooth initial data and different diffusion coefficients $\varepsilon \in \{1e{-3},1e{-4},1e{-5},1e{-6}\}$ (Modified Shu's Linear test with $T = 2\pi$, $\triangle t = 0.01$, $n_P = 500$, $m = 4$). The legend reports the discrete mass $M_h$ and energy $E_h$ at $t=0$ and $t=T$, which allows a direct comparison of mass variation and energy dissipation.}
    \label{fig:diffusion_effects}
\end{figure}

We evaluate the accuracy of our TBRBF method compared to the traditional Kansa method on various diffusion coefficients, \(\varepsilon \in \{1e{-3}, 1e{-4}, 1e{-5}, 1e{-6}\}\). As shown in \cref{fig:diffusion_effects}, when the diffusion coefficient \(\varepsilon\) drops below \(1e{-4}\), the final solution displays oscillations, especially near the corners of the trapezoidal region, leading to instability. In contrast, our proposed method yields stable numerical results. Moreover, the legend reports the discrete mass $M_h$ and energy $E_h$ at $t=0$ and $t=T$. By comparing these values, we observe that the mass variation is small, while the energy decreases over time. This provides a mass conservation and energy dissipation diagnostic for this non-smooth test.

\begin{figure}[htbp]
    \centering
    \begin{tabular}{c}
        \includegraphics[width=0.8\textwidth]{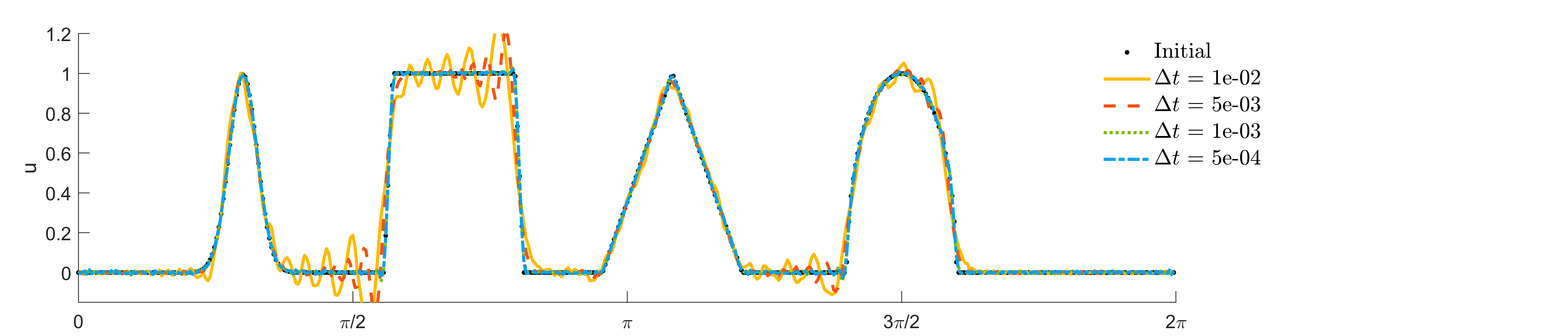}\\
        (a)  Traditional Kansa method solutions with different $\triangle t$\\
    (\small{$M_h(T) = M_h(0) \pm 4.3\times 10^{-3},  
E_h(T) = E_h(0) \pm 1.2\times 10^{-3}$})
       \\
         \vspace{1mm}
        \\
        \includegraphics[width=0.8\textwidth]{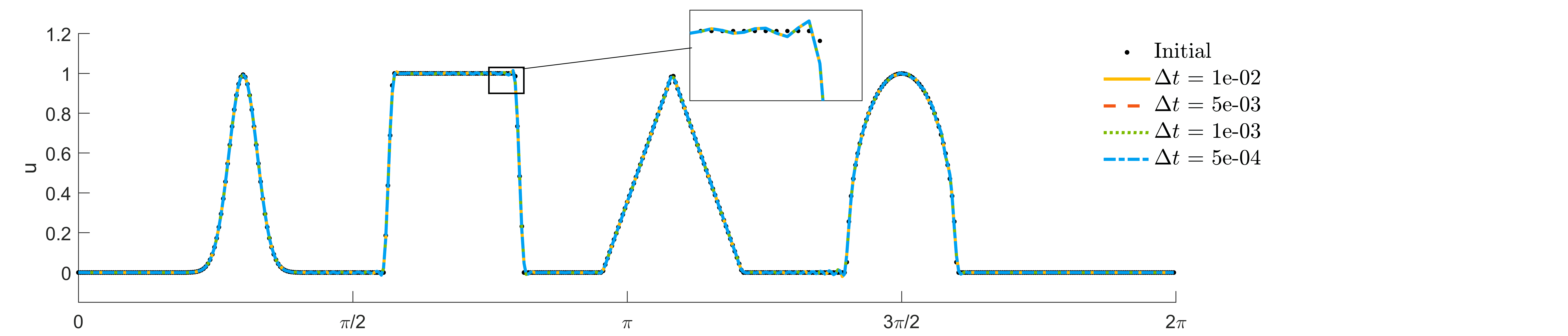}\\
        (b) TBRBF method solutions with different $\triangle t$\\
       (\small{$M_h(T) = M_h(0) \pm 4.4\times10^{-5}, 
E_h(T) = E_h(0) \pm 4.6\times10^{-4}$})
    \end{tabular}
    \caption{\textbf{(Example 2)} Performance comparison using  traditional Kansa method (top) and TBRBF method (bottom) to solve \eqref{eq:adv_diff} on the unit circle with non-smooth initial data and varying time steps 
             $\triangle t \in \{0.01,0.005,0.001,0.0005\}$ (Modified Shu's Linear test with $T = 2\pi$, $n_P = 500$, $\varepsilon = 1e{-6}$, $m = 4$).}
    \label{fig:dt_sensitivity}
\end{figure}

Semi-Lagrangian and characteristic treatments are not constrained by
the usual advective CFL condition in the same way as explicit Eulerian
discretizations. However, the stability of the fully discrete method
also depends on the trajectory integration, off-node evaluation, and
diffusion discretization. The following experiment therefore assesses
the observed stability of TBRBF for comparatively large time steps. For smaller values of diffusion coefficient \(\varepsilon = 10^{-6}\), we compared the stability of solutions over different time steps. 
As illustrated in \cref{fig:dt_sensitivity}, even with a relatively large time step of \(\triangle t = 10^{-2}\), our proposed method demonstrates good stability, while traditional methods require a time step of \(\triangle t \leq 10^{-3}\) to maintain stability. This characteristic opens up opportunities for simulating complex geometries in the future, enabling the acquisition of long time solutions while significantly reducing computation time and costs compared to traditional algorithms.





\subsection{Example 3: Toroidal dynamics with singularity data on the torus}
We extend the stability analysis from the unit circle to the surface experiment described in \cite{Ling-JSC22}, where the field \( u(x,t) \) evolves under a time-independent velocity field \( \mathbf{v} = (v_x, v_y, v_z) \)  on a torus with an inner radius of \( 1/3\) and an outer radius of 1. Details on the implicit surface representation and velocity expressions are provided in \cref{equations and velocities}. The velocity field specified in \cref{equations and velocities} advects particles along a (3,2) torus knot, returning them to their initial positions after a period \( T = 2\pi \). 

\begin{figure}[htbp]
    \centering
        \subfloat[Initial condition on the torus \qquad \\
         $M_h(0)\hspace{-1.2mm}=\hspace{-1.2mm}13.017042$, \hspace{-1.5mm} $E_h(0)\hspace{-1.2mm}=\hspace{-1.2mm}6.845094$ ]{
         \hspace{5mm} \includegraphics[width=0.22\textwidth] {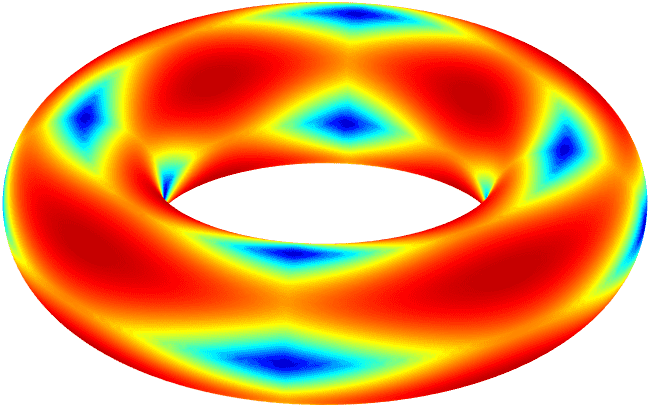} \hspace{5mm}
        \label{fig:torus_3D_init}
    }
       \hspace{20mm}
    \subfloat[Initial condition in the parameter space]{
        \hspace{10mm} \includegraphics[width=0.24\textwidth]{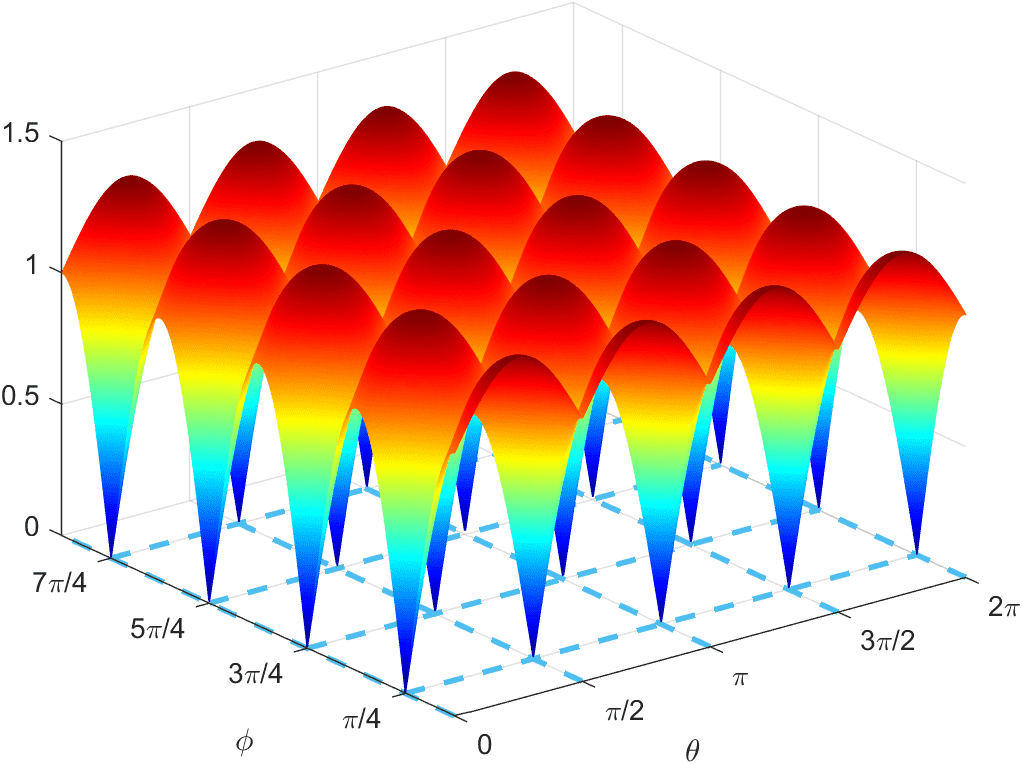} \hspace{10mm}
        \label{fig:param_space}
    }
    \caption{\textbf{(Example 3)} Visualization of the singular initial condition on the torus. Left: The initial condition $u(\theta, \phi)$ displayed on the torus surface. Right: The initial condition represented in the parameter space $(\theta, \phi)$, clearly highlighting the location and structure of the singularity (blue dashed). The values $M_h$ and $E_h$ shown below the panels are the discrete mass and energy at $t=0$.}
    \label{fig:initial_condition}
\end{figure}

In this case, as the independent variables of velocity are $\theta$ and $\phi$, we adopt a parametric representation and specify a singular initial condition
\begin{equation}
u(\theta, \phi) = |\sin(2\theta)| + |\cos(2\phi)| - 0.7|\sin(2\theta)\cos(2\phi)|. \label{eq:u_function}
\end{equation}
To aid visualization, \cref{fig:initial_condition} presents this initial condition: the left panel shows the condition mapped onto the torus, while the right one provides an intuitive view of the singularity in parameter space. For reference, we also report the discrete mass $M_h(0)$ and energy $E_h(0)$ of the initial data in the figure.

To emphasize the advection-dominated nature of the surface advection-diffusion equation \eqref{eq:adv_diff}, we set the diffusion coefficient to a small value \( \epsilon = 10^{-6} \), which can significantly impact the stability of the numerical solution. We employ an RK2-CN time-stepping scheme (with $T = 2\pi$ and $\triangle t = \frac{\pi}{500}$), while spatial discretization uses $n_P = 11600$ collocation points and the Mat\'{e}rn kernel with smoothness order $m = 6$.

\begin{figure}[tbhp]
   \centering
    \subfloat{%
        \includegraphics[width=0.6\textwidth]{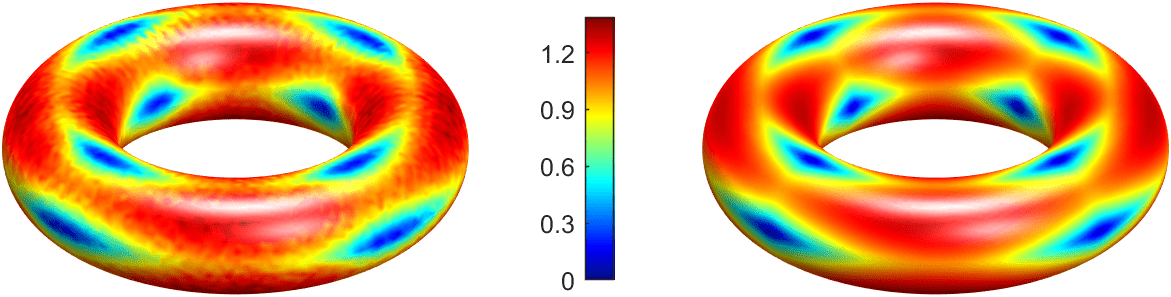}%
    }\\[2mm]
          {\small
\begin{tabular}{@{}ll@{}}
        (a) Solution of traditional Kansa method \qquad  \qquad & (b) Solution of TBRBF method\\
  {$M_h(T) = M_h(0)+8.71e-4$   }   & {$M_h(T) = M_h(0)+8.04e-4$}\\
  {$E_h(T) = E_h(0)+4.81e-3$  }   & { $E_h(T) = E_h(0)-3.57e-4$}
\end{tabular}
}\\
    \subfloat{%
        \includegraphics[width=0.6\textwidth]{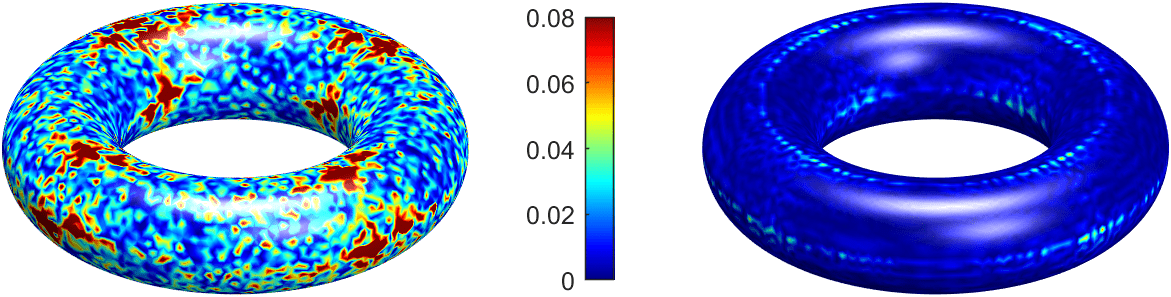}%
    }\\[2mm]
    {\small
      (c) Error of traditional Kansa method\qquad\qquad\:
      (d) Error of TBRBF method\quad
    }
    \caption{\textbf{(Example 3)} Numerical solutions and corresponding errors for solving \eqref{eq:adv_diff} on a torus with a (3,2) torus knot velocity. The top-left panel shows the solution obtained using the traditional Kansa method, while the top-right panel shows the solution obtained using the TBRBF method. The bottom-left and bottom-right panels display the corresponding  errors with respect to the reference solution for the traditional Kansa method and the proposed method, respectively. Discretization parameters used are $(T, \triangle t, n_P, m) = (2\pi, \pi/500, 11600, 6)$. The values $M_h$ and $E_h$ shown below the solution panels are the discrete mass and energy at $t=T$.}
    \label{fig:dynamic_evolution}
\end{figure}

Given the small value of $\varepsilon$ and the fact that $T=2\pi$ represents a complete rotation, we use the initial condition as a reference solution. As shown in \cref{fig:dynamic_evolution}, the conventional Kansa method performs poorly, producing significant inaccuracies due to the singularity and the small value of $\varepsilon$. In contrast, our proposed method maintains the integrity of the solution's shape and achieves high accuracy, even with the singular initial condition. We also report the discrete mass $M_h$ and energy $E_h$ at $t=T$ in \cref{fig:dynamic_evolution}. We note that in this example $\varepsilon=10^{-6}$ is very small and the initial data has corner-type singularities. 
Both effects make the problem more sensitive to discretization errors. 
As a result, the discrete energy $E_h$ is not guaranteed to decrease monotonically in this setting. 
This behavior is more visible for the traditional Kansa method. 
For the TBRBF method the energy change is much smaller and the mass variation remains small.

\subsection{Example 4: Turing patterns model with advection terms}
We now illustrate the effectiveness of the proposed TBRBF method by applying it to reaction–advection–diffusion systems on a variety of general surfaces, which serve as benchmark problems for surface PDEs \cite{Shankar-SISC20}. Following the numerical experiments in \cite{Liu-CMA25}, we consider the two-species Turing system with advection velocity $\mathbf{v}$ as defined in \cref{equations and velocities}, posed on the sphere, torus, Bretzel2, and CPD surfaces (each implicitly defined by algebraic equations in \cref{equations and velocities})
\begin{equation}
\begin{aligned}
\partial_t u + \nabla_{\mathcal{M}} \cdot (u\mathbf{v}) &= \varepsilon_u \Delta_{\mathcal{M}} u + a u (1 - k_1 w^2) + w (1 - k_2 u), \\
\partial_t w + \nabla_{\mathcal{M}} \cdot (w\mathbf{v}) &= \varepsilon_w \Delta_{\mathcal{M}} w + b w\left( 1 + \frac{a k_1}{b}u w\right) + u (c + k_2 w).
\end{aligned}
\label{reaction-advection-diffusion system}
\end{equation}
The model parameters depend on the specific surface but satisfy the relations $\varepsilon_u = 0.516\varepsilon_w$, $a = -c = 0.899$, and $b = -0.91$. The remaining parameters, including those for the spots and stripes on each surface, are provided in \cref{tab:parameters1}. The initial condition is generated by assigning uniformly random values in the range $[-0.5, 0.5]$ within a narrow strip along the surface's equator, with zeros elsewhere.
\begin{table}[htb]
   \centering
   \caption{\textbf{(Example 4)} The parameters of the Turing system used to simulate the spots and stripes with $\varepsilon_u = 0.516\varepsilon_w$, $b = -0.91$, $a = -c = 0.899$.  }
   \label{tab:parameters1}
   \begin{tabular}{cccccc} 
   \toprule
   Surface & Pattern & $\varepsilon_w$ & $k_1$ & $k_2$ & Final time \\
   \midrule
   \multirow{2}{*}{Sphere}        & Spots   & \(6.50 e{-3}\) & 0.02 & 0.2  & 600 \\
                                  & Stripes & \(8.87e {-4}\) & 3.5  & 0    & 6000 \\
    \midrule
   \multirow{2}{*}{Torus}         & Spots   & \(3.50 e{-3}\) & 0.02 & 0.2  & 1200 \\
                                  & Stripes & \(8.87 e{-4}\) & 3.5  & 0    & 6000 \\
    \midrule
   \multirow{2}{*}{Bretzel2}  & Spots   & \(2.50 e{-3}\) & 0.02 & 0.2  & 1000 \\
                                  & Stripes & \(1.80 e{-4}\) & 3.5  & 0    & 6000 \\
    \midrule
   \multirow{2}{*}{CPD surface}           & Spots   & \(2.50e{-3}\) & 0.02 & 0.2  & 600 \\
                                  & Stripes & \(3.50 e{-4}\) & 3.5  & 0    & 8000 \\      
   \bottomrule
   \end{tabular}
\end{table}

We apply the proposed TBRBF method to simulate the Turing patterns model, utilizing the RK2-CN scheme with a time step of $\triangle t = 0.5$ for temporal discretization and the Mat\'{e}rn kernel with smoothness order $m = 4$ for spatial discretization. Since all surfaces under consideration are defined implicitly, the normal vectors $\mathbf{\hat{N}}$ can be calculated exactly. For the closest point computations required by the surface-restricted RK method, we directly use MATLAB’s built-in nonlinear optimization solver, configured with the sequential quadratic programming (SQP) algorithm.

\begin{figure}[tbhp]
    \centering
    {\includegraphics[width=.18\linewidth]{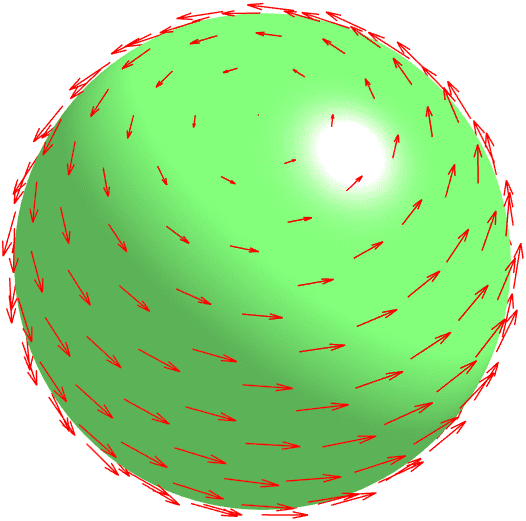}
     \hspace{11mm}
    \includegraphics[width=.18\linewidth]{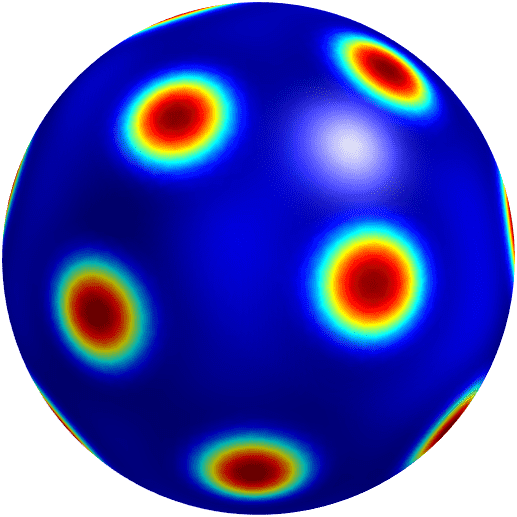}
    \hspace{11mm}
     \includegraphics[width=.18\linewidth]{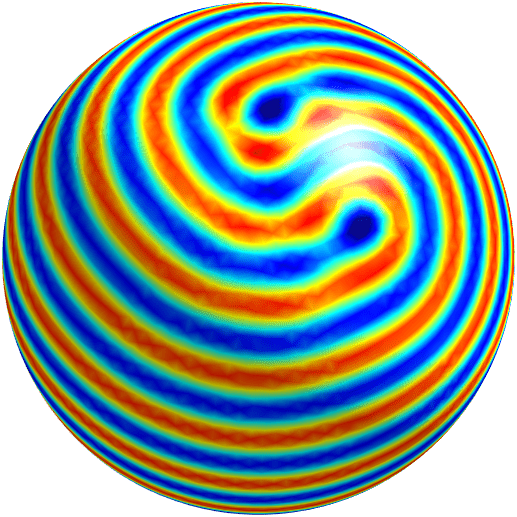}}
      \vspace{10mm}
     \\
    {\includegraphics[width=.2\linewidth]{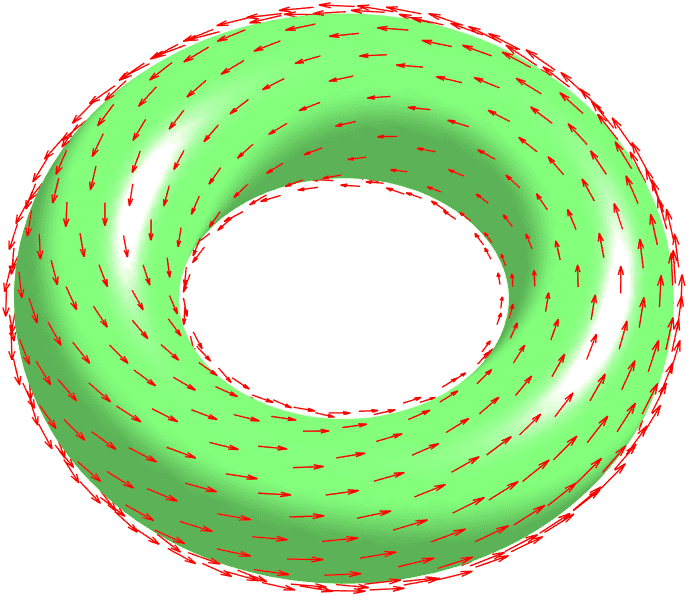}
    \hspace{7mm}
    \includegraphics[width=.2\linewidth]{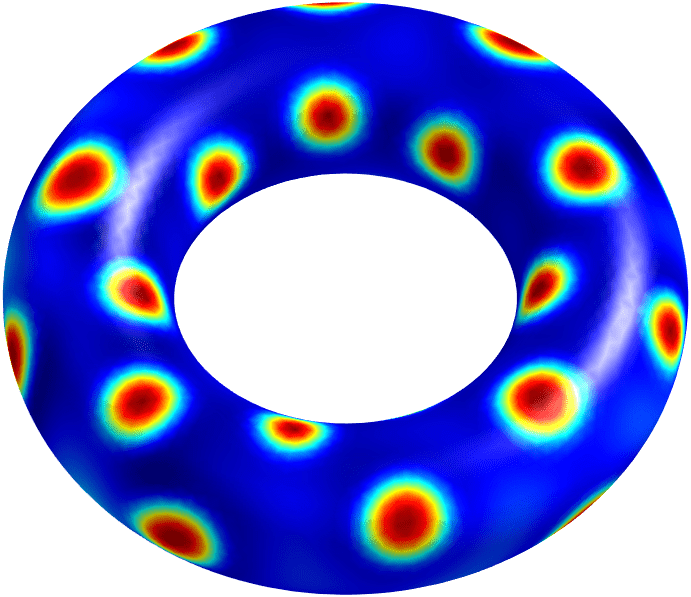}
    \hspace{7mm}
     \includegraphics[width=.2\linewidth]{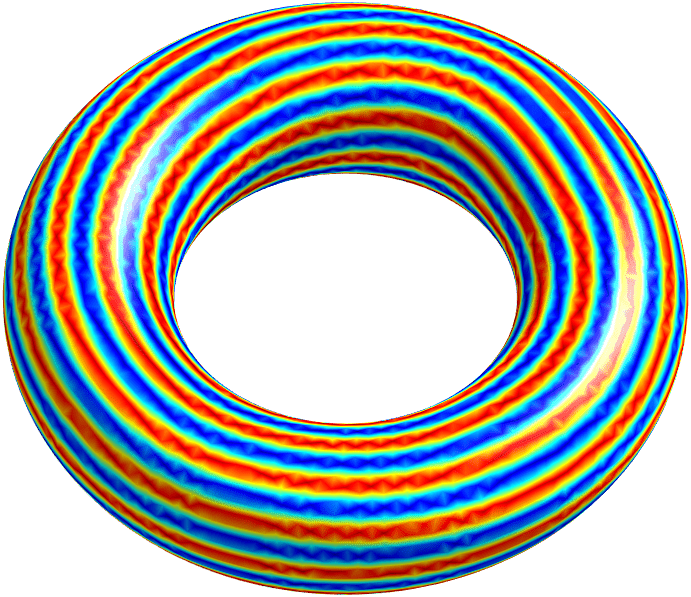}}
      \vspace{10mm}
     \\
   {\includegraphics[width=.22\linewidth]{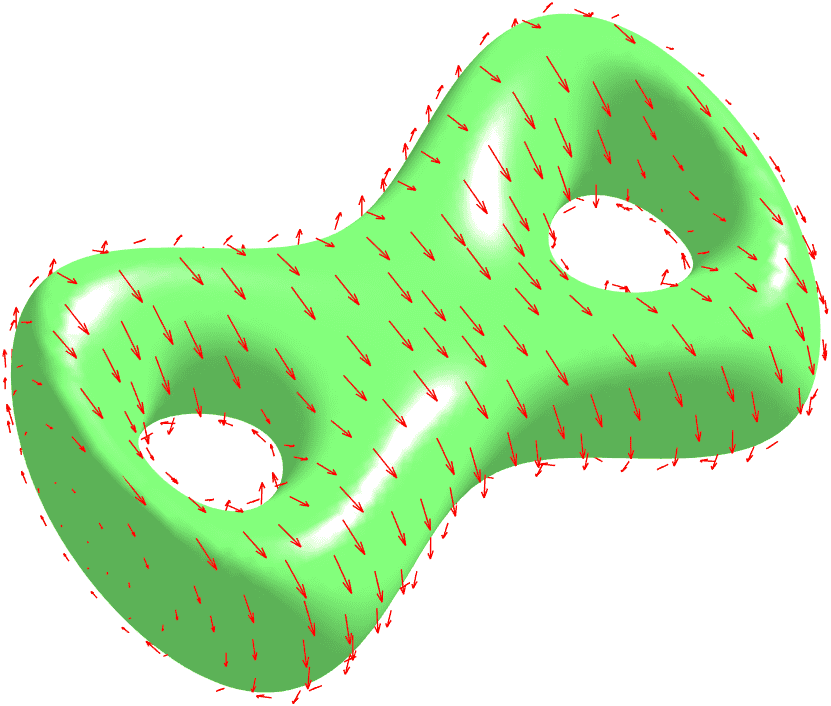}
   \hspace{5mm}
   \includegraphics[width=.22\linewidth]{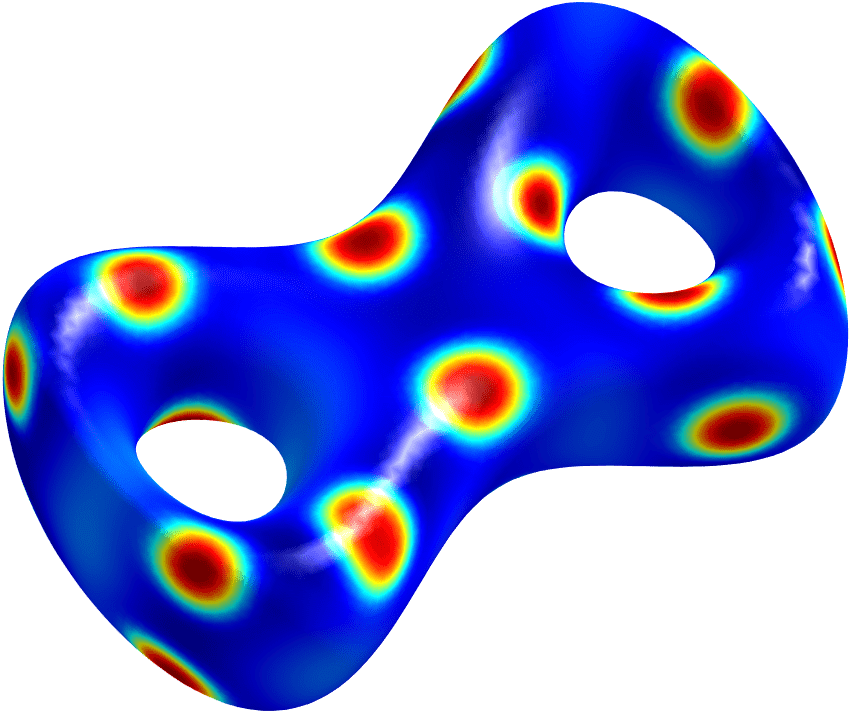}
   \hspace{5mm}
     \includegraphics[width=.22\linewidth]{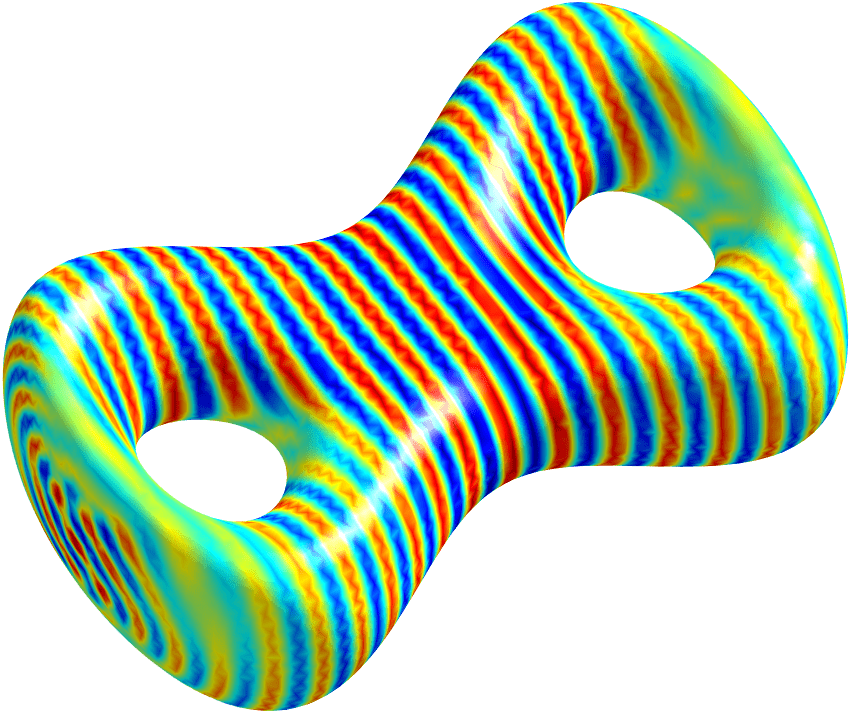}} 
      \vspace{10mm}
     \\  
   {\includegraphics[width=.2\linewidth]{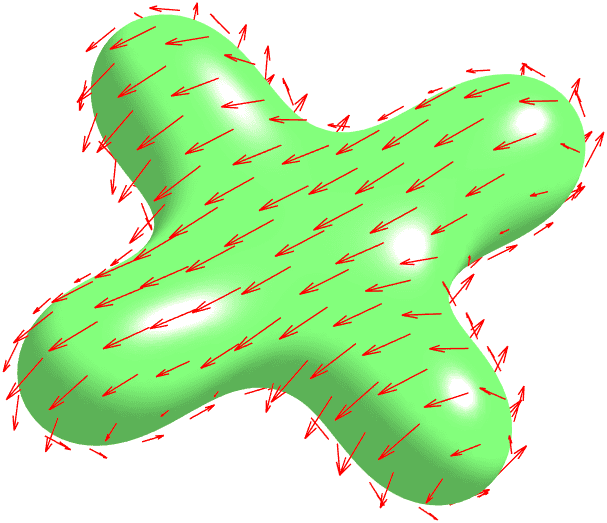}
   \hspace{5mm}
   \includegraphics[width=.2\linewidth]{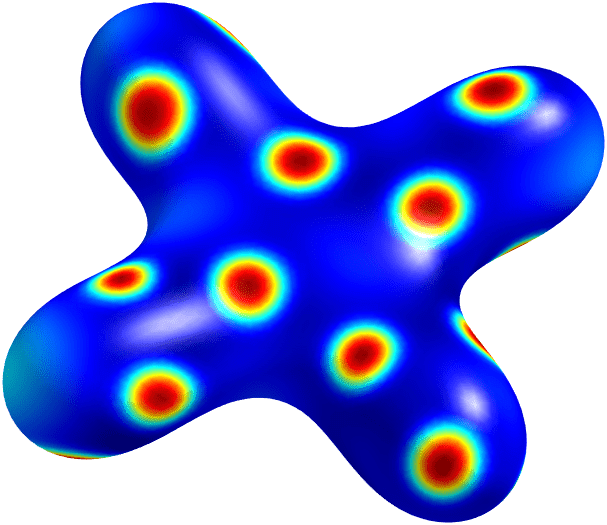}
   \hspace{10mm}
     \includegraphics[width=.2\linewidth]{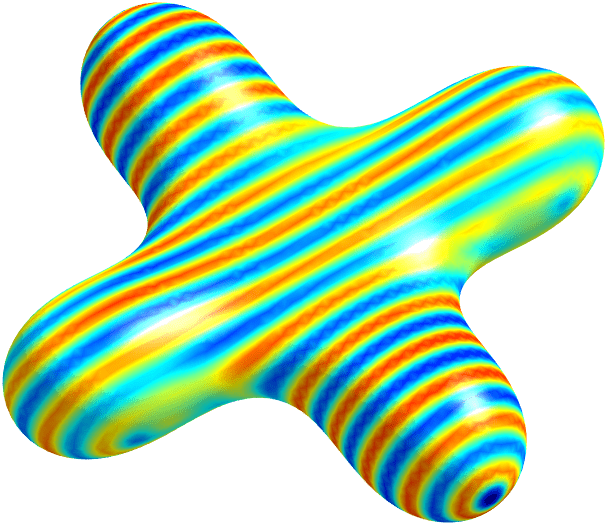}}
    \caption{\textbf{(Example 4)} Spots and Stripes on four distinct implicit surfaces: Sphere ($n_P = 4482$), Torus ($n_P = 5312$), Bretzel2 ($n_P = 7216$) and CPD surface ($n_P = 6256$). In each row, the left panel shows the velocity vector field, while the middle and right panels display the spots and stripes, respectively.}
    \label{fig:pattern_evolution}
\end{figure}

The simulation results are presented in \cref{fig:pattern_evolution}, where lighter and darker shades indicate higher and lower concentrations, respectively. The first column of \cref{fig:pattern_evolution} visualizes the advection velocity field, allowing observation of the transport dynamics for both the Spots and Stripes systems. As shown in \cref{fig:pattern_evolution}, our method successfully generates dynamic Spots that are advected across the surface by the flow. For the Stripes system, while the localized effects of the velocity field on the Bretzel2 and CPD surfaces prevent the development of sharply defined boundaries, the flow nonetheless induces a clear alignment of the patterns along its direction.

\subsection{Example 5: Point Clouds Extension for Turing Patterns}
To further demonstrate the versatility of our approach, we apply the TBRBF method, using the same discretization schemes and parameter settings as subsection 5.4 Example 4, to more general point cloud surfaces that lack explicit implicit or parametric representations. Specifically, we test on the Stanford Bunny and cow models. The surface reaction–advection–diffusion system \eqref{reaction-advection-diffusion system} is used, with model parameters provided in \cref{tab:parameters2}.
\begin{table}[htb]
   \centering
   \caption{\textbf{(Example 5)} The parameters of the Turing system used to simulate the spots and stripes on the Stanford Bunny and Cow with $\varepsilon_u = 0.516\varepsilon_w$, $b = -0.91$, $a = -c = 0.899$. }
   \label{tab:parameters2}
   \begin{tabular}{cccccc} 
   \toprule
   Surface & Pattern & $\varepsilon_w$ & $k_1$ & $k_2$ & Final time \\
   \midrule
   \multirow{2}{*}{Stanford Bunny \&  Cow}           & Spots   & \(1.50 e{-3}\) & 0.02 & 0.2  & 600 \\
                                  & Stripes & \(8.87e-4\) & 3.0  & 0    & 6000 \\
                                  
   \bottomrule
   \end{tabular}
\end{table}

Unlike implicit surfaces, exact normal vectors are not directly available for point clouds. Therefore, we estimate the normal vector at each point using localized level-set interpolation. For each point in the cloud, the estimated local normal is first aligned with the $z$-axis through a suitable rotation, and the local surface is reconstructed using an intrinsic method. The closest point is subsequently refined using Newton’s method. Further implementation details can be found in \cite{Ling-JSC22}.

\begin{figure}[htbp]
    \centering
    \includegraphics[width=0.23\textwidth]{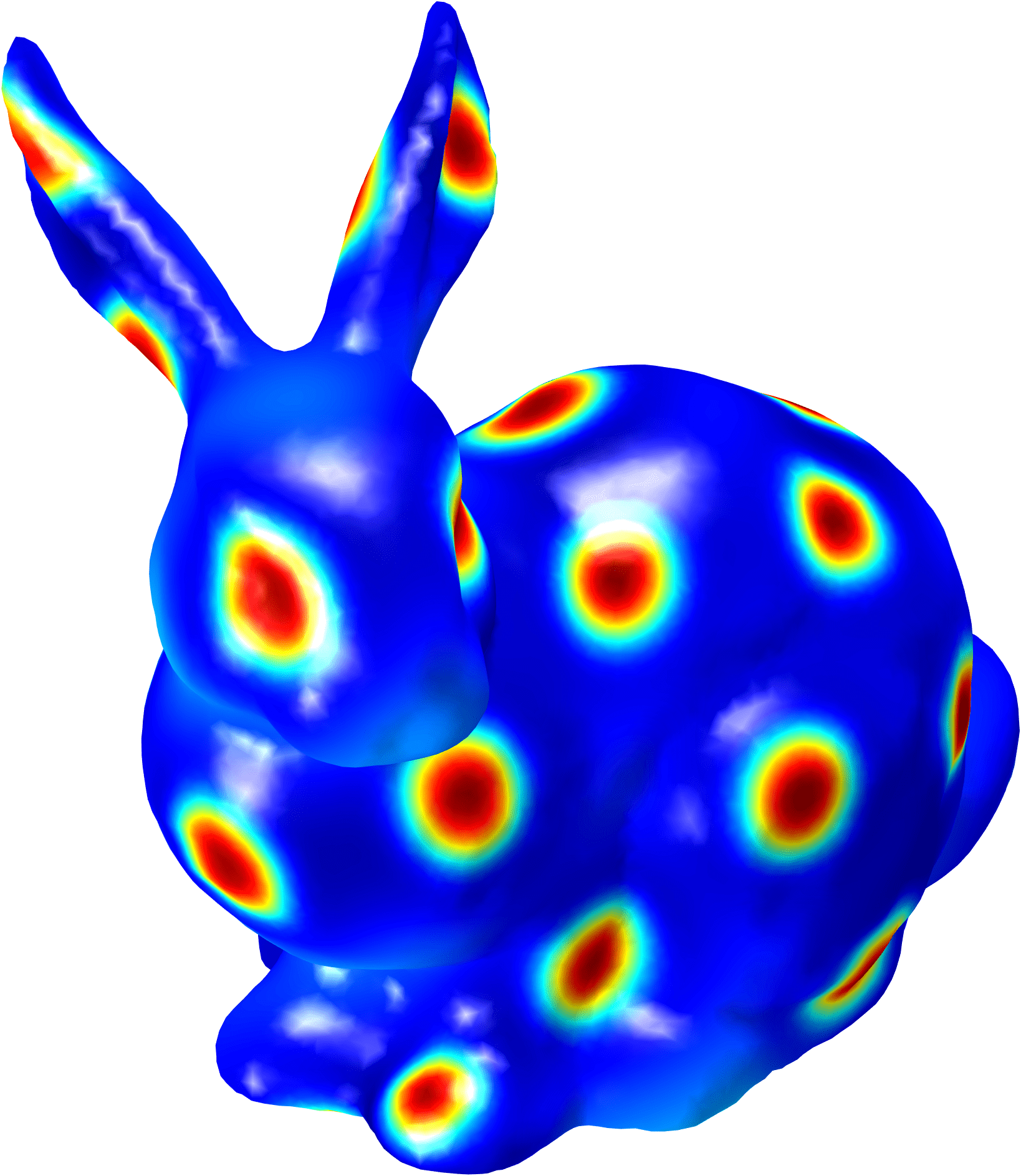}
    \hspace{15mm}
    \includegraphics[width=0.23\textwidth]{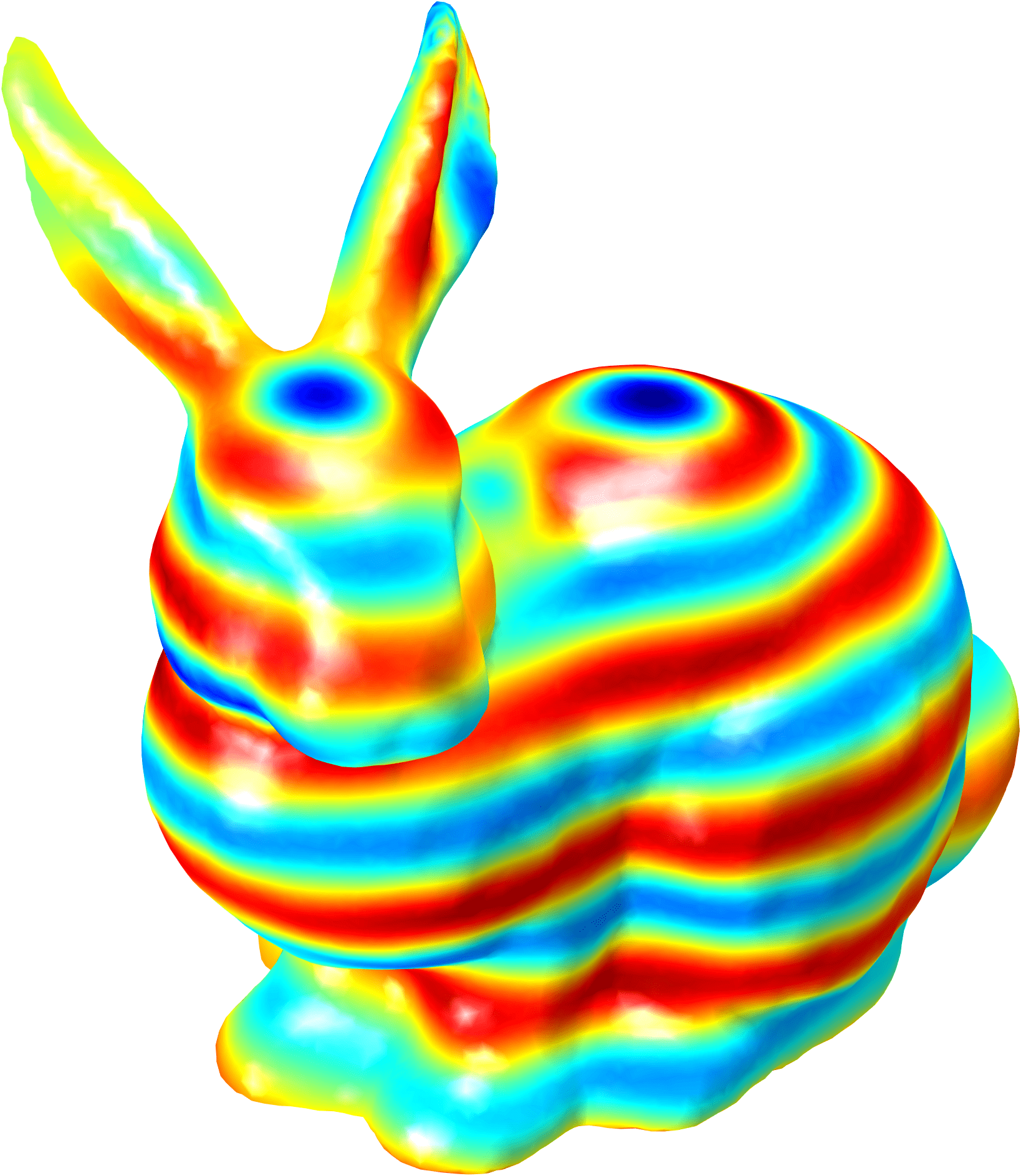}
    \vspace{5mm}
    \\
    \includegraphics[width=0.23\textwidth]{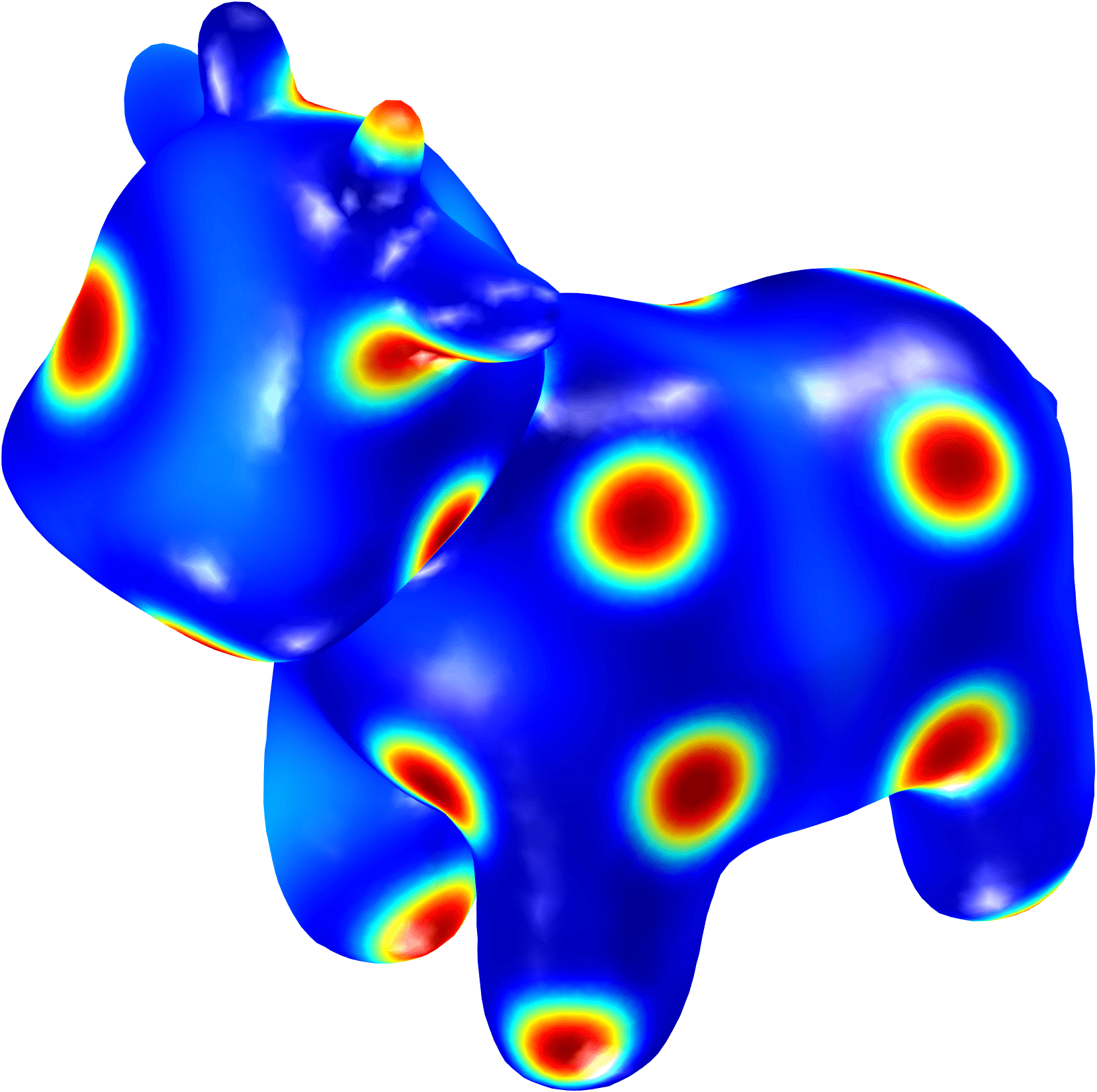}
    \hspace{15mm}
    \includegraphics[width=0.23\textwidth]{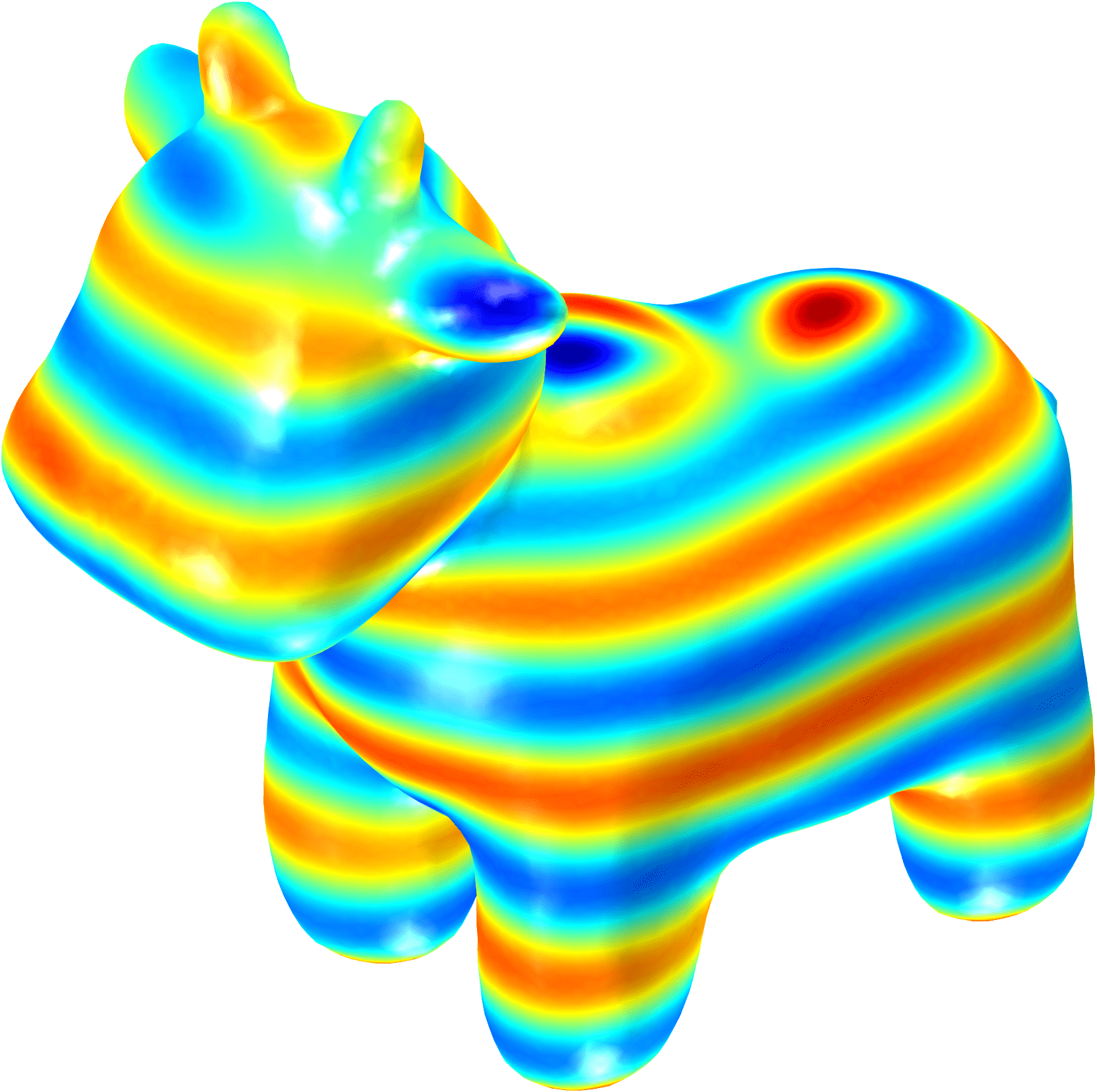}
    \caption{\textbf{(Example 5)} Simulation of Turing patterns on point cloud surfaces using the proposed TBRBF method. Results are shown for Stanford Bunny ($n_P = 15685$) and Cow ($n_P = 14358$). }
    \label{fig:point_cloud}
\end{figure}

With these procedures, we obtain the simulation results shown in \cref{fig:point_cloud}. Similar to the patterns observed in \cref{fig:pattern_evolution}, our method successfully captures the characteristic spot and stripe patterns of the Turing system on complex point cloud surfaces. The dynamic evolution and organization of these patterns are well preserved, demonstrating the robustness and adaptability of our approach even in the absence of explicit surface representations.

\section{Conclusion}
In this paper, we present a TBRBF method for solving the surface advection-diffusion equation on a closed, smooth, compact codimension-1 manifold. By integrating Kansa-type RBF collocation with a trajectory-based framework, our approach directly processes the Laplace–Beltrami operator on the manifold while tracking the evolution of the solution along the advection velocity trajectories. This method reduces the need for repeated narrow-band interpolation associated with fixed-time CAN enforcement. The surface formulation avoids a separate narrow-band CAN
re-extension stage. Its computational cost is instead associated with
trajectory projection, off-node RBF evaluation, and the solution of
the collocation system. We rigorously established the equivalence between the original surface advection-diffusion PDE and the OSC surface system by embedding the problem into a narrow band. Numerical experiments further validate the enhanced stability and accuracy of our proposed method. Overall, this work provides a robust computational strategy for handling surface PDEs while establishing a solid theoretical foundation for surface trajectory-based formulations.

Additionally, this article offers insights into using time-continuous embedded methods to tackle time-dependent problems, such as heat equations, advection diffusion reaction problems, and Navier-Stokes equations. This approach provides a way to reduce repeated CAN enforcement and interpolation back to the manifold in trajectory-based surface computations. Looking ahead, we aim to extend this method to solve Navier-Stokes equations and thin shell problems. However, it is important to note that the proposed method currently struggles with shock and oscillation problems, and our future work will focus on addressing these challenges.

\appendix
\section{Differential-geometric identities in an orthonormal chart}\label{app:diffgeo}
In this appendix we collect standard surface differential-geometric definitions used in Section~2.
We work in a local chart $\Gamma=\Phi^{-1}$ with coordinates $(\theta_1,\theta_2)$.
The metric matrix $G_{\Gamma}$ of the first fundamental form is defined by
\begin{equation}\label{G1stFF}
    G_{\Gamma} = \begin{bmatrix}
        E  & F \\ F  & G
    \end{bmatrix}
    = \begin{bmatrix}
        \partial_{\theta_1} \Gamma \cdot \partial_{\theta_1} \Gamma  & \partial_{\theta_1} \Gamma \cdot \partial_{\theta_2} \Gamma  \\ 
        \partial_{\theta_1} \Gamma \cdot \partial_{\theta_2} \Gamma  & \partial_{\theta_2} \Gamma \cdot \partial_{\theta_2} \Gamma
    \end{bmatrix}.
\end{equation}
In our advection-aligned orthonormal coordinates \eqref{Dtheta}, we have $G_{\Gamma}=\mathbf I$.
To characterize the curvature of $\mathcal{M}\in \mathbb{R}^3$, we can use the inner product of the second-order partial derivatives of $\Gamma(\theta_1,\theta_2)$ with the unit normal vector $\mathbf{\hat N}$ to obtain following curvature matrix $H_{\Gamma}$ of the second fundamental form
\begin{equation}\label{H2stFF}
    H_{\Gamma} = \begin{bmatrix}
        L  & M \\ M  & N
    \end{bmatrix}
    = \begin{bmatrix}
        \partial_{\theta_1\theta_1} \Gamma \cdot \mathbf{\hat N}  & \partial_{\theta_1\theta_2} \Gamma \cdot \mathbf{\hat N}  \\ 
        \partial_{\theta_1\theta_2} \Gamma \cdot \mathbf{\hat N}  & \partial_{\theta_2\theta_2} \Gamma \cdot \mathbf{\hat N}
    \end{bmatrix}.
\end{equation}
Let $S_{\Gamma}$ denote the shape operator, which is defined as
\begin{equation}\label{ShapeO}
    S_{\Gamma} = G_{\Gamma}^{-1} H_{\Gamma}.
\end{equation}
In the orthonormal chart, $G_\Gamma=\mathbf I$, hence $S_\Gamma=H_\Gamma$.
Finally, to express the derivatives of the normal vector in the direction of any tangent vector along our advection-aligned coordinate directions, we can derive the following Weingarten equation
\begin{equation}\label{Weingarten-app}
 \Big[\partial_{\theta_1} \mathbf{\hat N}, \,\partial_{\theta_2} \mathbf{\hat N}\Big]  
 = - \Big[  \partial_{\theta_1}  \Gamma, \, \partial_{\theta_2}  \Gamma \Big] S_{\Gamma}
 = - \Big[  \partial_{\theta_1}  \Gamma, \, \partial_{\theta_2}  \Gamma \Big]\begin{bmatrix} L & M \\ M & N \end{bmatrix}.
\end{equation}
Equation \eqref{Weingarten-app} reduces to \eqref{Weingarten} in Section~2 since $G_\Gamma=\mathbf I$.

\section{Implicit algebraic equations and advection velocities for all tested manifolds}\label{equations and velocities}

The implicit algebraic equations for general surfaces are defined
\small{ \begin{align*}
&\bullet\ \text{Sphere}:\quad \left\{(x, y,z) \in \mathbb{R}^3 \mid x^2 + y^2 +z^2 = 1\right\}. \\
&\bullet\ \text{Torus}: \quad \left\{(x, y, z) \in \mathbb{R}^3 \mid (x^2 + y^2 + z^2 + 1^2 - (1/3)^2)^2 - 4(x^2 + y^2) = 0 \right\}. \\
&\bullet\ \text{Bretzel2}: \quad \left\{(x, y, z) \in \mathbb{R}^3 \mid (x^2(1 - x^2) + y^2)^2 + 1/2z^2  - 1/40(x^2 + y^2 + z^2) = 1/40 \right\}. \\
&\bullet\ \text{CPD surface}: \quad \left\{(x, y, z) \in \mathbb{R}^3 \mid \sqrt{(x - 1)^2 + y^2 + z^2}\,\sqrt{(x + 1)^2 + y^2 + z^2} \right.\\
&\qquad\qquad\qquad\qquad\qquad\qquad \left.\, \sqrt{x^2 + (y - 1)^2 + z^2}\,\sqrt{x^2 + (y + 1)^2 + z^2} = 1.1 \right\}.
\end{align*}}

The corresponding advection velocities are given as follows
\small{\begin{align*}
&\bullet\ \text{Sphere:} \quad \mathbf{v} = 0.01\, \mathbf{\hat N} \times (0,0,-1), &&\\
&\bullet\ \text{Torus:} \quad 
\begin{cases} 
v_x = \rho_2 \cos(3\theta) - 3\rho_1 \sin(3\theta),\\
v_y = \rho_2 \sin(3\theta) + 3\rho_1 \cos(3\theta),\\
v_z = 2r \cos(2\phi),
\end{cases} &&\\
&\quad \hspace{14mm} \text{where } \rho = \sqrt{x^2+y^2}, \theta = \tfrac{1}{3}\tan^{-1}\Bigl(\frac{y}{x}\Bigr), \phi = \tfrac{1}{2}\tan^{-1}\Bigl(\frac{z}{\rho-R}\Bigr),&&\\
&\quad \hspace{23mm} R = 1, ~ r = \tfrac{1}{3}, ~ \rho_1 = R + r \cos(2\phi), ~ \rho_2 = -2r \sin(2\phi), &&\\
&\bullet\ \text{Bretzel2:}\quad \mathbf{v} = 0.01\, \mathbf{\hat N} \times (-1,0,0), &&\\
&\bullet\ \text{CPD surface:}\quad \mathbf{v} = 0.01\, \mathbf{\hat N} \times (0,1,0), &&\\
&\bullet\ \text{Point clouds:} \quad \mathbf{v} = 0.01\, \mathbf{\hat N} \times (0,0,1). &&
\end{align*}}

\bibliographystyle{siamplain}
\bibliography{references}
\end{document}